\documentclass[12pt,twoside]{article}
\usepackage{mathrsfs}
\usepackage{amsfonts}
\usepackage{multirow}
\usepackage[dvips]{graphicx}
\usepackage{mathrsfs}
\usepackage{amsfonts}
\usepackage{graphicx}
\usepackage{subfigure}
\usepackage{bbm}
\usepackage{mathrsfs}
\usepackage{dsfont}
\usepackage{mathrsfs}
\usepackage{amsmath}
\usepackage{amssymb}
\usepackage{amsfonts}

\usepackage{amsmath,amssymb}
\usepackage{mathrsfs,color}
\newenvironment{proof}[1][Proof]{\noindent\textbf{#1.} }{\ \rule{0.5em}{0.5em}}
\newtheorem{theorem}{\quad Theorem}[section]
\newtheorem{definition}{\quad Definition}[section]

\newtheorem{example}{\quad Example}[section]
\newtheorem{lemma}{\quad Lemma}[section]
\numberwithin{figure}{section} \numberwithin{equation}{section}
\makeatletter 
\newcommand{\bm}[1]{\mbox{\boldmath{$#1$}}}
\usepackage{amsmath}
\begin{document}
\date{}
\author{Hengfei Ding$^{1}$\thanks{E-mail:
dinghf05@163.com},
\;\;Changpin Li$^{2}$\thanks{E-mail:
lcp@shu.edu.cn},\;\;Qian Yi$^{2}$\\
 \small \textit{ 1. School of Mathematics and Statistics, Tianshui
Normal University},
\\ \small\textit{ Tianshui 741001, China}\\
 \small \textit{2. Department of Mathematics, Shanghai University, Shanghai 200444, China}\\
 }
\title{A new second-order midpoint approximation formula for Riemann-Liouville derivative:
 algorithm and its application\thanks{The work was partially supported by the National
Natural Science Foundation of China under Grant Nos. 11372170 and 11561060,
 the Scientific Research Program
 for Young Teachers of Tianshui Normal University under Grant No. TSA1405, and Tianshui
Normal University Key Construction Subject Project (Big data processing in dynamic image).
}}\maketitle
 \hrulefill
\begin{abstract}
Compared to the the classical first-order Gr\"{u}nwald-Letnikov formula at time $t_{k+1} (\textmd{or}\; t_{k})$, we firstly
 propose a second-order numerical approximate scheme for discretizing the
Riemann-Liouvile derivative
at time $t_{k+\frac{1}{2}}$, which is very suitable for constructing the Crank-Niclson technique
 applied to the time-fractional differential equations. The established formula has the following form
$$
\begin{array}{lll}
\displaystyle \,_{\mathrm{RL}}{{{\mathrm{D}}}}_{0,t}^{\alpha}u\left(t\right)\left|\right._{t=t_{k+\frac{1}{2}}}= \tau^{-\alpha}\sum\limits_{\ell=0}^{k}
\varpi_{\ell}^{(\alpha)}u\left(t_k-\ell\tau\right)
+\mathcal{O}(\tau^2),\,\,k=0,1,\ldots, \alpha\in(0,1),
\end{array}
$$
where the coefficients $\varpi_{\ell}^{(\alpha)}$ $(\ell=0,1,\ldots,k)$
can be determined via the following generating function
 $$
\begin{array}{lll}
\displaystyle G(z)=\left(\frac{3\alpha+1}{2\alpha}-\frac{2\alpha+1}{\alpha}z+\frac{\alpha+1}{2\alpha}z^2\right)^{\alpha},\;|z|<1.
\end{array}
$$
Applying this formula to the time fractional Cable equations
with Riemann-liouville derivative in one or two space dimensions.
Then the high-order compact finite difference schemes are obtained.
The solvability, stability and convergence with orders
$\mathcal{O}(\tau^2+h^4)$ and $\mathcal{O}(\tau^2+h_x^4+h_y^4)$ are shown,
where $\tau$ is the temporal stepsize and $h$, $h_x$, $h_y$ are the spatial stepsizes, respectively. Finally,
numerical experiments are provided to support the theoretical analysis.
\\ \vspace{0.2 cm}
\vspace{0.2 cm}\\
\end{abstract}
\hrulefill

\section{Introduction}

 In recent years, a great deal of attention has been focused on fractional differential equations
 due to well describing many physical
processes and phenomenons \cite{BDA,BWM,MK1,MK2,MRGS,SK}.
Very limited analytical methods, such as the Fourier transform method, the Laplace transform method
and the Green function method are used to solve the very special fractional differential equations.
 So to seek numerical methods is the center task for studies of fractional differential equations
 \cite{A,DL1,DL2,ER,JLPZ,LD,PCSCJ,SLLA,WV,WY,ZLTA}. In the history of numerical methods for fractional differential equations,
Liu et al. \cite{LAT}, and Meerschaert and Tadjeran \cite{MT} are the first ones that developed the finite difference methods
for fractional partial differential equations.
The Galerkin finite element methods for fractional partial differential
equations is proposed by Ervin and Roop, for the stationary space fractional partial differential equations with two-sided
Riemann-Liouville derivatives. They first presented a rigorous analysis of the well-posedness
of the weak formulation in the framework of fractional Sobolev spaces \cite{ER}.

Generally speaking, one of key issues of approximating fractional differential equations is how to numerically discretice
the fractional derivatives. Although there have existed some studies on numerical approximations of factional
integrals and fractional derivatives, high-order scheme for time fractional derivatives have not been throughly solved.
This paper aims to construct new and effective the second-order mid-point approximate formula for time Riemann-Liouville
derivative. Then the established scheme is applied to time fractional Cable equations in one and two space dimensions.
From bibliography available, there have existed numerical studies for the fractional Cable equations. For example,
Langlans et al. \cite{LHW} developed two implicit finite difference schemes with convergence orders
 $\mathcal{O}(\tau + h^2)$ and $\mathcal{O}(\tau^2 + h^2)$.
Hu and Zhang proposed two implicit compact difference schemes, where the first scheme was proved to be stable and convergent with
order $O(\tau + h^4)$ by the energy method \cite{HZ}.
In \cite{MY}, Quintana-Murillo and Yuste constructed an explicit numerical scheme for
fractional Cable equation which includes two temporal Riemann-Liouville derivatives, where they showed the
stability and convergence conditions by using the Von Neumann method.
Zhuang et al. \cite{ZLTA} considered the one-dimensional time fractional Cable
equation by using the Galerkin finite element method, in which
the proposed method was based on a semi-discrete finite difference approximation
in time and Galerkin finite element method in space.
The spectral
method for fractional Cable equation was discussed by Lin et al. \cite{LLX}, where the detailed theoretical analysis was provided.
As far as we know, the computational efficiency for time fractional Cable equation is not high yet. Besides, the high-dimensional
time fractional Cable equations seen not to be studied. Here, we study the fractional cable equation in two space dimensions
where the fractional derivative is approximated by the derived method in this paper. The unconditional stability and convergence of
the established numerical algorithms are presented by the energy method.

The reminder of the paper is constructed as follows. In Section 2, we establish a new second-order
approximation formula for Riemann-Liouville derivatives. Then
two high-order finite difference schemes for the fractional Cable equations in one and two space dimensions
are proposed in Sections 3 and 4, respectively. Numerical experiments are displayed in
Section 5, where are in line with the theoretical analysis. Remarks and conclusions are included in the last section.

\section{Second-order scheme for Riemann-Liouville derivative}

 In the section, we propose a new second-order approximation formula for computing Riemann-Liouville derivatives.

\begin{definition} \cite{P,SKM}
Let $-\infty \leq a < t < b \leq \infty$. For given $u\in L_1([a,b])$, its
 Riemann-Liouville derivative
of order $\alpha>0$ with lower limit $a$ is defined as follows
$$
\begin{array}{lll}
\displaystyle \,_{\mathrm{RL}}{{{\mathrm{D}}}}_{a,t}^{\alpha}u(t) =
\frac{1}{\Gamma(n-\alpha)}\frac{\mathrm{d}^n}{\mathrm{d}t^n}
 \int_{0}^{t}\frac{u(s)}{(t-s)^{\alpha-n+1}}\mathrm{d}s,\;\; n-1<\alpha \leq n\in \mathds{Z}^{+}.
\end{array}
$$
\end{definition}

\begin{lemma} \cite{P,SKM} Suppose $\alpha>0$, $u(t)\in L_1(\mathds{R})$. Then the Fourier transform of
 $\alpha$-th Riemann-Liouville derivative of $u(t)$ is
$$
\begin{array}{lll}
\displaystyle \mathcal{F}\left\{{\,_{\mathrm{RL}}\mathrm{D}}_{-\infty,t}^{\alpha}u(t)\right\}(\omega) =(-\mathrm{i}\omega)^{\alpha}\widehat{u}(\omega),
\end{array}
$$
where $\widehat{u}(\omega)=\int_{\mathds{R}}\mathrm{e}^{\mathrm{i}\omega t}u(t)\mathrm{d}t$ denotes the Fourier transform of $u(t)$.
\end{lemma}

Now, we start to develop the second-order numerical approximation formula.
\begin{theorem}\label{th:3.1} Denote $$\mathscr{C}^{n+\alpha}(\mathds{R})=\left\{u|u\in L_1(\mathds{R}),
\;\int_{\mathds{R}}\left(1+|\omega|\right)^{n+\alpha}|\widehat{u}(\omega)|\mathrm{d}\omega<\infty\right\}.$$
 Define the following difference operator
$$
\begin{array}{lll}
\displaystyle \delta_{-\infty,t}^{\alpha}u\left(t\right)
=\tau^{-\alpha}\sum\limits_{\ell=0}^{\infty}
\varpi_{\ell}^{(\alpha)}u\left(t-\ell\tau\right),
\end{array}\eqno(1)
$$
where $\tau$ is temporal stepsize. If $u\in\mathscr{C}^{2+\alpha}(\mathds{R})$,
then one has
$$
\begin{array}{lll}
\displaystyle \,_{\mathrm{RL}}{{{\mathrm{D}}}}_{-\infty,t}^{\alpha}u\left(t+\frac{1}{2}\tau\right)
=\mathcal{\delta}_{-\infty,t}^{\alpha}u\left(t\right)+\mathcal{O}(\tau^2)
\end{array}\eqno(2)
$$
as $\tau\rightarrow0$.

Here $\displaystyle \varpi_{\ell}^{(\alpha)}\;\left(
\ell=0,1,\ldots,\right)$ are the expansion coefficients of $G(z)$, that is,
$$
\begin{array}{lll}
\displaystyle G(z)=
\sum\limits_{\ell=0}^{\infty}\varpi_{\ell}^{(\alpha)}z^\ell,\;\;|z|<1.
\end{array}\eqno(3)
$$
where
 $$G(z)=\left(\frac{3\alpha+1}{2\alpha}-\frac{2\alpha+1}{\alpha}z+\frac{\alpha+1}{2\alpha}z^2\right)^{\alpha}.$$
\end{theorem}

\begin{proof}
Taking the Fourier transform on both sides of equation (1) then combining with equation (3), one has
$$
\begin{array}{lll}
\displaystyle \mathcal{F}\left\{\delta_{-\infty,t}^{\alpha}u\left(t\right)\right\}(\omega)
&=&\displaystyle\tau^{-\alpha}\sum\limits_{\ell=0}^{\infty}
\varpi_{\ell}^{(\alpha)}\mathrm{e}^{\mathrm{i}\omega\ell\tau}\hat{u}(\omega)\vspace{0.2 cm}\\
&=&\displaystyle\tau^{-\alpha}\hat{u}(\omega)
\sum\limits_{\ell=0}^{\infty}
\varpi_{\ell}^{(\alpha)}\mathrm{e}^{\mathrm{i}\omega\ell\tau}\vspace{0.2 cm}\\
&=&\displaystyle (-\mathrm{i}\omega)^{\alpha}\mathrm{e}^{-\mathrm{i}\frac{\omega\tau}{2}}S(\mathrm{i}\omega\tau)\hat{u}(\omega),
\end{array}
$$
where $S(z)=(-z)^{-\alpha}\mathrm{e}^{\frac{z}{2}}G(\mathrm{e}^{z})$.

Note that
$$
\begin{array}{lll}
\displaystyle G(\mathrm{e}^{z})&=&\displaystyle\left(\frac{3\alpha+1}{2\alpha}-\frac{2\alpha+1}{\alpha}
\mathrm{e}^{z}+\frac{\alpha+1}{2\alpha}\mathrm{e}^{2z}\right)^{\alpha}\vspace{0.2 cm}\\
&=&\displaystyle\left(-z\right)^{\alpha}\left(1-\frac{1}{2\alpha}z-\frac{2\alpha+3}{6\alpha}z^2+\cdots\right)^{\alpha}\vspace{0.2 cm}\\
&=&\displaystyle(-z)^{-\alpha}\left(1-\frac{1}{2}z-\frac{8\alpha^2+9\alpha+3}{24\alpha}z^2+\cdots\right)
\end{array}
$$
and
$$
\begin{array}{lll}
\displaystyle S(z)&=&\displaystyle(-z)^{-\alpha}\mathrm{e}^{\frac{z}{2}}G(\mathrm{e}^{z})=\left(1+\frac{z}{2}+\frac{z^2}{8}+\cdots\right)
\left(1-\frac{1}{2}z-\frac{8\alpha^2+9\alpha+3}{24\alpha}z^2+\cdots\right)\vspace{0.2 cm}\\
&=&\displaystyle1-\frac{8\alpha^2+12\alpha+3}{24\alpha}z^2+\mathcal{O}(|z|^3).
\end{array}
$$
So, there exists a constant $c_1>0$ such that
$ \left|1-S(\mathrm{i}\omega\tau)\right|\leq c_1|\omega\tau|^2$.

Then
$$
\begin{array}{lll}
\displaystyle \mathcal{F}\left\{\delta_{-\infty,t}^{\alpha}u\left(t\right)\right\}(\omega)
&=&\displaystyle(-\mathrm{i}\omega)^{\alpha}\mathrm{e}^{-\mathrm{i}\frac{\omega\tau}{2}}\hat{u}(\omega)+
 (-\mathrm{i}\omega)^{\alpha}\mathrm{e}^{-\mathrm{i}\frac{\omega\tau}{2}}
 \left(S(\mathrm{i}\omega\tau)-1\right)\hat{u}(\omega)\vspace{0.2 cm}\\&=&\displaystyle
 \mathcal{F}\left\{\,_{\textmd{RL}}{{{\mathrm{D}}}}_{-\infty,t}^{\alpha}
 u\left(t+\frac{1}{2}\tau\right)\right\}(\omega)+\hat{\varphi}(\tau,\omega)
\end{array}
$$
by using Lemma 2.1.
Hence, one has
$$
\begin{array}{lll}
\displaystyle |{\varphi}(\tau,\omega)|=\left|\frac{1}{2\pi \mathrm{i}}
\int_{-\infty}^{\infty}\textmd{e}^{-\mathrm{i}\omega\tau}\hat{\varphi}(\tau,\omega)\textmd{d}\omega\right|\leq
\frac{c_1\tau^2}{2\pi}\int_{-\infty}^{\infty}(1+|\omega|)^{2+\alpha}|\hat{u}(\omega)|\textmd{d}\omega=\mathcal{O}(\tau^2).
\end{array}
$$
All this ends the proof.
\end{proof}

{\bf Remark:} If $u(t)$ is suitably smooth and has compact support for $0<t<T$, then Riemann-Liouville
derivative $\,_{\textmd{RL}}{\mathrm{D}}_{-\infty,t}^{\alpha}u(t)$ coincides with
$\,_{\textmd{RL}}{\mathrm{D}}_{0,t}^{\alpha}u(t)$. Denote
$$
\begin{array}{lll}
\displaystyle \delta_{0,t}^{\alpha}u\left(t\right)
=\tau^{-\alpha}\sum\limits_{\ell=0}^{\left[\frac{t}{\tau}\right]}
\varpi_{\ell}^{(\alpha)}u\left(t-\ell\tau\right).
\end{array}
$$
Then
the corresponding numerical approximation formula (2) is reduced to
$$
\begin{array}{lll}
\displaystyle \,_{\textmd{RL}}{\mathrm{D}}_{0,t}^{\alpha}u\left(t\right)\left|\right._{t=t_{k+\frac{1}{2}}}
=\delta_{0,t}^{\alpha}u\left(t_k\right)+\mathcal{O}(\tau^2),\;\;k=0,1,\ldots,
\end{array}
$$
that is,
$$
\begin{array}{lll}
\displaystyle \,_{\mathrm{RL}}{{{\mathrm{D}}}}_{0,t}^{\alpha}u\left(t\right)
\left|\right._{t=t_{k+\frac{1}{2}}}= \tau^{-\alpha}\sum\limits_{\ell=0}^{k}
\varpi_{\ell}^{(\alpha)}u\left(t_k-\ell\tau\right)
+\mathcal{O}(\tau^2),\,\,k=0,1,\ldots
\end{array}\eqno(4)
$$
the coefficients $\varpi_\ell^{(\alpha)}$ $(\ell=0,1,\ldots)$ in equation (3) can be expressed as follows,
$$
\begin{array}{lll}
 \displaystyle \varpi_{\ell}^{(\alpha)}=
\left(\frac{3\alpha+1}{2\alpha}\right)^{\alpha}\sum\limits_{m=0}^{\ell}
\left(\frac{\alpha+1}{3\alpha+1}\right)^m g_{m}^{(\alpha)}
g_{\ell-m}^{(\alpha)},\;\;\ell=0,1,\ldots,
\end{array}\eqno(5)
$$
where $$\displaystyle g_{m}^{(\alpha)}=(-1)^m
\left(\alpha\atop
m\right)=(-1)^m\frac{\Gamma(1+\alpha)}{\Gamma(m+1)\Gamma(1+\alpha-m)}
,\;m=0,1,\ldots,\ell.$$

For convenience, take the place of $\alpha$ in $\varpi_{\ell}^{(\alpha)}$ by $1-\alpha$ which can not cause confusion, where $\alpha\in(0,1)$.
It is easy to get the following theorem.

\begin{theorem}  The coefficients $\varpi_{\ell}^{(1-\alpha)},\;\ell=0,1,\ldots.$ can be computed
recursively by the formulas.
$$\left\{
\begin{array}{lll}
\displaystyle
\varpi_{0}^{(1-\alpha)}=\left(\frac{4-3\alpha}{2(1-\alpha)}\right)^{1-\alpha};\vspace{0.2 cm}\\
\displaystyle \varpi_{1}^{(1-\alpha)}=\frac{2(\alpha-1)(3-2\alpha)}{4-3\alpha}\varpi_{0}^{(1-\alpha)};\vspace{0.2 cm}\\
\displaystyle \varpi_{\ell}^{(1-\alpha)}=\frac{1}{(4-3\alpha)\ell}\left[2(\ell+\alpha-2)(3-2\alpha)\varpi_{\ell-1}^{(1-\alpha)}
+(2-\alpha)(4-2\alpha-\ell)\varpi_{\ell-2}^{(1-\alpha)}\right],\;\ell\geq2.
\end{array}\right.
$$
\end{theorem}

\begin{theorem}  The coefficients $\varpi_{\ell}^{(1-\alpha)} $ are nonpositive if $\ell\geq5$ where $\alpha\in(0,1)$, that is,
$\varpi_{\ell}^{(1-\alpha)}\leq0$ if $\ell\geq5$, where $\alpha\in(0,1)$.
\end{theorem}

\begin{proof}
For $0<\alpha<1$ and $\ell\geq5$, one has
$$
\begin{array}{lll}
 \displaystyle \varpi_{\ell}^{(1-\alpha)}&=&\displaystyle
\left(\frac{4-3\alpha}{2(1-\alpha)}\right)^{1-\alpha}\sum\limits_{m=0}^{\ell}
\left(\frac{2-\alpha}{4-3\alpha}\right)^m g_{m}^{(1-\alpha)}
g_{\ell-m}^{(1-\alpha)}\vspace{0.2 cm}\\
&=&\displaystyle\left\{\left[1+\left(\frac{2-\alpha}{4-3\alpha}\right)^\ell\right]
g_{0}^{(1-\alpha)}
g_{\ell}^{(1-\alpha)}+\left[\frac{2-\alpha}{4-3\alpha}+\left(\frac{2-\alpha}{4-3\alpha}\right)^{\ell-1}\right]g_{1}^{(1-\alpha)}
g_{\ell-1}^{(1-\alpha)}\right.\vspace{0.2 cm}\\&&\displaystyle\left.+\sum\limits_{m=2}^{\ell-2}
\left(\frac{2-\alpha}{4-3\alpha}\right)^m g_{m}^{(1-\alpha)}
g_{\ell-m}^{(1-\alpha)}\right\}\left(\frac{4-3\alpha}{2(1-\alpha)}\right)^{1-\alpha}\vspace{0.2 cm}\\
&\leq&\displaystyle\left\{\left[1+\left(\frac{2-\alpha}{4-3\alpha}\right)^\ell\right]
g_{0}^{(1-\alpha)}
g_{\ell}^{(1-\alpha)}+\left[\frac{2-\alpha}{4-3\alpha}+\left(\frac{2-\alpha}{4-3\alpha}\right)^{\ell-1}\right]g_{1}^{(1-\alpha)}
g_{\ell-1}^{(1-\alpha)}\right.\vspace{0.2 cm}\\&&\displaystyle\left.+ g_{2}^{(1-\alpha)}
g_{\ell-2}^{(1-\alpha)}\sum\limits_{m=2}^{\infty}
\left(\frac{2-\alpha}{4-3\alpha}\right)^m \right\}\left(\frac{4-3\alpha}{2(1-\alpha)}\right)^{1-\alpha}
\vspace{0.2 cm}\\
&=&\displaystyle\left\{\left[1+\left(\frac{2-\alpha}{4-3\alpha}\right)^\ell\right]
\left(1-\frac{2-\alpha}{\ell}\right)
+\left[\frac{2-\alpha}{4-3\alpha}+\left(\frac{2-\alpha}{4-3\alpha}\right)^{\ell-1}\right](\alpha-1)\right.\vspace{0.2 cm}\\&&\displaystyle\left.-
\frac{\alpha(\ell-1)(2-\alpha)^2}{4(4-3\alpha)(\alpha+\ell-3)} \right\}\left(\frac{4-3\alpha}{2(1-\alpha)}\right)^{1-\alpha}g_{\ell-1}^{(1-\alpha)}
\vspace{0.2 cm}\\
&\leq&\displaystyle\left\{\left[1+\left(\frac{2-\alpha}{4-3\alpha}\right)^\ell\right]
\left(1-\frac{2-\alpha}{\ell}\right)
+\left[\frac{2-\alpha}{4-3\alpha}+\left(\frac{2-\alpha}{4-3\alpha}\right)^{\ell-1}\right](\alpha-1)\right.\vspace{0.2 cm}\\
\end{array}
$$
$$
\begin{array}{lll}
 \displaystyle
&=:&\displaystyle\left\{P(\alpha,\ell)+ Q(\alpha,\ell)
\left(\frac{2-\alpha}{4-3\alpha}\right)^{\ell}\right\}g_{\ell-1}^{(1-\alpha)}.
\end{array}
$$
Here
$$
\begin{array}{lll}
\displaystyle P(\alpha,\ell)=\left[
\frac{\ell-2+\alpha}{\ell}+\frac{(\alpha-1)(2-\alpha)}{4-3\alpha}-\frac{\alpha(\ell-1)(2-\alpha)^2}{4(4-3\alpha)(\alpha+\ell-3)}
-\frac{3(1-\alpha)}{2}\left(\frac{2-\alpha}{4-3\alpha}\right)^{4}
\right],
\end{array}
$$
and
$$
\begin{array}{lll}
\displaystyle Q(\alpha,\ell)=\left[\frac{\ell-2+\alpha}{\ell}+\frac{3(1-\alpha)}{2}+\frac{(\alpha-1)(4-3\alpha)}{2-\alpha}\right].
\end{array}
$$

It is somewhat tedious but easy to
check that $P(\alpha,\ell)$ and $Q(\alpha,\ell)$ are increasing with respect to $\ell$. Hence,
$$
\begin{array}{lll}
&\displaystyle P_{\min}(\alpha,\ell)=\displaystyle P(\alpha,5)\vspace{0.2 cm}\\
&=\displaystyle
\frac{177\alpha^6-969\alpha^5+3990\alpha^4-12000\alpha^3+20000\alpha^2-16336\alpha+5152}{10(2+\alpha)(4-3\alpha)^4}>0,
\end{array}
$$
and
$$
\begin{array}{lll}
\displaystyle Q_{\min}(\alpha,\ell)=Q(\alpha,5)=
\frac{-17\alpha^2+23\alpha+2}{10(2-\alpha)}>0.
\end{array}
$$
Noticing that $g_{\ell-1}^{(1-\alpha)}\leq0$ for $0<\alpha<1$ and $\ell\geq 1$ gives
$$
\begin{array}{lll}
\displaystyle
\varpi_{\ell}^{(1-\alpha)}&\leq&\displaystyle
\left\{P(\alpha,\ell)+ Q(\alpha,\ell)
\left(\frac{2-\alpha}{4-3\alpha}\right)^{\ell}\right\}g_{\ell-1}^{(1-\alpha)}
\vspace{0.2 cm}\\
&\leq&\displaystyle
\left\{P_{\min}(\alpha,\ell)+ Q_{\min}(\alpha,\ell)
\left(\frac{2-\alpha}{4-3\alpha}\right)^{\ell}\right\}g_{\ell-1}^{(1-\alpha)}\leq0.
\end{array}
$$
All this completes the proof.
\end{proof}

\begin{theorem} The coefficient $\varpi_{\ell}^{(1-\alpha)} $ are increasing with
 respect to $\ell\geq7$, where $\alpha\in(0,1)$, that is,
$\varpi_{\ell}^{(1-\alpha)}\geq\varpi_{\ell-1}^{(1-\alpha)}$ if $\ell\geq 7$, where $\alpha\in(0,1)$.
\end{theorem}

\begin{proof}
Here, we use mathematical induction to prove this theorem. Let
$$
\begin{array}{ll}
\displaystyle
S_{\ell}=\varpi_{\ell}^{(1-\alpha)}-\varpi_{\ell-1}^{(1-\alpha)},\;\;\ell\geq7.
\end{array}
$$

 From equation (5), one easily knows that $S_{7}\geq0$.
 Now suppose that the conclusion holds for $\ell= k>7$, that is
 $$
\begin{array}{ll}
\displaystyle
\varpi_{k}^{(1-\alpha)}\geq\varpi_{k-1}^{(1-\alpha)},\;\; k>7.
\end{array}
$$
 Then for $\ell=k+1$, according to Theorems 2.3 and 2.4, one gets
 $$
\begin{array}{lll}
\displaystyle
S_{k+1}&=&\displaystyle\varpi_{k+1}^{(1-\alpha)}-\varpi_{k}^{(1-\alpha)}\vspace{0.2 cm}\\
&=&\displaystyle\frac{1}{(k+1)(4-3\alpha)}\left[2(k+\alpha-1)(3-2\alpha)\varpi_{k}^{(1-\alpha)}\right.\vspace{0.2 cm}\\&&\displaystyle\left.
+(2-\alpha)(3-2\alpha-k)\varpi_{k-1}^{(1-\alpha)}-(k+1)(4-3\alpha)\varpi_{k}^{(1-\alpha)}\right]
\vspace{0.2 cm}\\
&\geq&\displaystyle\frac{1}{(k+1)(4-3\alpha)}\left[2(k+\alpha-1)(3-2\alpha)
+(2-\alpha)(3-2\alpha-k)\right.\vspace{0.2 cm}\\&&\displaystyle\left.-(k+1)(4-3\alpha)\right]\varpi_{k}^{(1-\alpha)}
\vspace{0.2 cm}\\
&=&\displaystyle\frac{2(1-\alpha)(\alpha-2)}{(k+1)(4-3\alpha)}\varpi_{k}^{(1-\alpha)}\geq0.
\end{array}
$$
The proof is thus completed.
\end{proof}

\section{Application to the fractional Cable equation in one space dimension}

 In this section, we study the following one-dimensional Cable equation
$$
\begin{array}{rrr}
\displaystyle
\frac{\partial u(x,t)}{\partial t}=K_1\,_{\mathrm{RL}}{{{\mathrm{D}}}}_{0,t}^{1-\alpha_1}
\frac{\partial^2u(x,t) }{\partial x^2}
-K_2 \,_{\mathrm{RL}}{{{\mathrm{D}}}}_{0,t}^{1-\alpha_2} u(x,t) +f(x,t),\vspace{0.2 cm}\\\;(x,t)\in(0,L)\times(0,T],
\end{array}\eqno(6)
$$
with initial condition
$$
\begin{array}{ll}
\displaystyle
u(x,0)=0,\;\;x\in[0,L],
\end{array}\eqno(7)
$$
and boundary value conditions
$$
\begin{array}{ll}
\displaystyle
u(0,t)=\varphi_1(t),\;u(L,t)=\varphi_2(t),\;\;t\in(0,T],
\end{array}\eqno(8)
$$
where $0<\alpha_1,\alpha_2<1$, $K_1$ and $K_2$ are two constants,
$f(x, t)$, $\varphi_1(t)$ and $\varphi_2(t)$ are suitably smooth functions.

\subsection{Development of numerical algorithm}

Firstly, denote $x_j = jh$, $t_k = k\tau$,
$\Omega_h = \{x_j| 0\leq j\leq M\},
\Omega_\tau= \{t_k|0\leq k\leq N\},$ and $\Omega_{\tau,h}=\Omega_\tau\times\Omega_h$,
where $h =L/M,
\tau=T/N$
 are the uniform spatial and temporal mesh sizes respectively, and $M, N$ are
two positive integers. Let
$$
\begin{array}{ll}
\displaystyle
\mathcal{S}_h=\{\bm{u}|\bm{u}=(u_0,u_1,\ldots,u_M),\;u_0=u_M=0\}
\end{array}
$$
be defined on $\Omega_h$.

In addition, define the following first- and second-order difference operators as,
$$
\begin{array}{ll}
\displaystyle
\delta_x u_{j-\frac{1}{2}}=\frac{1}{h}(u_j-u_{j-1}),\;\;\delta_x^2 u_{j}=\frac{1}{h}\left(\delta_x u_{j+\frac{1}{2}}-\delta_x u_{j-\frac{1}{2}}\right),
\end{array}
$$
and fourth-order compact difference operator $\mathscr{L}$ as,
$$\mathscr{L}u_j^k=\left\{
\begin{array}{ll}
\left(I+\frac{h^2}{12}\delta_x^2\right) u_{j}^k=\frac{1}{12}(u_{j+1}^k+10u_j^k+u_{j-1}^k),\;\;1\leq j\leq M-1,\vspace{0.2 cm}\\
u_j^k,\;\;j=0,M,
\end{array}\right.
$$
where $I$ is the unit operator.

Now we turn to derive an effective finite difference scheme for solving equation (6),
together with initial and boundary value conditions (7) and (8).
 Consider equation (6) at point $\left(x_j,t_{k+\frac{1}{2}}\right)$
$$
\begin{array}{ll}
\displaystyle
\frac{\partial u\left(x_j,t_{k+\frac{1}{2}}\right)}{\partial t}=&\displaystyle K_1
\,_{\mathrm{RL}}{{{\mathrm{D}}}}_{0,t}^{1-\alpha_1}\left(
\frac{\partial^2u\left(x_j,t_{k+\frac{1}{2}}\right) }{\partial x^2}\right)\vspace{0.2 cm}\\
&\displaystyle- K_2\,_{\mathrm{RL}}{{{\mathrm{D}}}}_{0,t}^{1-\alpha_2} u\left(x_j,t_{k+\frac{1}{2}}\right)
 +f\left(x_j,t_{k+\frac{1}{2}}\right).
\end{array}\eqno(9)
$$

Applying the second-order central difference formula
$$
\begin{array}{ll}
\displaystyle
\frac{\partial u\left(x_j,t_{k+\frac{1}{2}}\right)}{\partial t}=
\frac{u\left(x_j,t_{k+1}\right)-u\left(x_j,t_{k}\right)}{\tau}+\mathcal{O}(\tau^2),
\end{array}
$$
and second-order approximation formula (4) to the above equation (9), one gets
$$
\begin{array}{ll}
\displaystyle
\frac{u\left(x_j,t_{k+1}\right)-u\left(x_j,t_{k}\right)}{\tau}=&\displaystyle K_1\delta_{0,t}^{1-\alpha_1}
\frac{\partial^2u\left(x_j,t_{k}\right) }{\partial x^2}
-K_2 \delta_{0,t}^{1-\alpha_2}u\left(x_j,t_{k}\right)\vspace{0.2 cm}\\&
 +f\left(x_j,t_{k+\frac{1}{2}}\right)+\mathcal{O}(\tau^2).
\end{array}\eqno(10)
$$
Acting the operator $\mathscr{L}$ on both sides of (10) and noticing
$$
\begin{array}{ll}
\displaystyle
\mathscr{L}\frac{\partial^2u\left(x_j,t_{k}\right) }{\partial x^2}=\delta_x^2u\left(x_j,t_{k}\right) +\mathcal{O}(h^4),
\end{array}
$$
one has
$$
\begin{array}{ll}
\displaystyle
\mathscr{L}(u\left(x_j,t_{k+1}\right)-u\left(x_j,t_{k}\right))=&\displaystyle\tau K_1\delta_{0,t}^{1-\alpha_1}\delta_x^2u\left(x_j,t_{k}\right)
-\tau K_2\delta_{0,t}^{1-\alpha_2}\mathscr{L}u\left(x_j,t_{k}\right)\vspace{0.2 cm}\\&\displaystyle
+\tau\mathscr{L}f\left(x_j,t_{k+\frac{1}{2}}\right)+\tau R_j^k,
\end{array}
$$
in which there exists a positive
constant $C_1$  such that
$|R_j^k|\leq C_1(\tau^2+h^4)$.

Omit the local truncation error $R_j^k$ and denote the numerical
solution of $u(x_j,t_k)\in\Omega_{\tau,h}$ by $u_j^k$. One can establish the following
 high-order compact difference scheme for equation (6), together with (7) and (8),
$$
\begin{array}{rrr}
\displaystyle
\mathscr{L}\left(u_j^{k+1}-u_j^{k}\right)=\tau K_1\delta_{0,t}^{1-\alpha_1}\delta_x^2u_j^{k}
-\tau K_2\delta_{0,t}^{1-\alpha_2}\mathscr{L}u_j^{k}  +\tau\mathscr{L}f_j^{k+\frac{1}{2}},\vspace{0.2 cm}\\\;\;1\leq j\leq M-1,\;1\leq k\leq N-1,\vspace{0.2 cm}\\
\end{array}\eqno(11)
$$
$$
\begin{array}{lll}
\displaystyle
\displaystyle
u_j^0=0,\;\;0\leq j\leq M,
\end{array}\eqno(12)
$$
$$
\begin{array}{lll}
\displaystyle
u_0^k=\varphi_1(t_k),\;u_M^k=\varphi_2(t_k),\;\;1\leq k\leq N.
\end{array}\eqno(13)
$$

\subsection{Solvability, stability and convergence analysis}

For arbitrary vectors $\bm{u},\bm{v}\in\mathcal{S}_h$, we introduce the following inner products and the corresponding norms,
$$
\begin{array}{lll}
\displaystyle
(\bm{u},\bm{v})=h\sum_{j=1}^{M-1}u_jv_j,\;\;||\bm{u}||=\sqrt{(\bm{u},
\bm{u})},\;\;(\delta_x^2\bm{u},\bm{v})=h\sum_{j=1}^{M-1}(\delta_x^2u_j)v_j,\\
\displaystyle
(\delta_x\bm{u},\delta_x\bm{v})=h\sum_{j=1}^{M}(\delta_x u_j)(\delta_x v_j),
\;\;||\delta_x\bm{u}||=\sqrt{(\delta_x\bm{u},\delta_x\bm{u})}.
\end{array}
$$
It is easy know $(\delta_x^2\bm{u},\bm{v})=-(\delta_x\bm{u},\delta_x\bm{v})$.

Next several lemmas are listed which will be used later on.
\begin{lemma}\cite{WG} \label{Lem.3.1} If $a\geq0$, $b\geq0$, then the following Young's inequality holds,
$$
\begin{array}{lll}
\displaystyle
ab\leq\frac{(\epsilon a)^p}{p}+\frac{\left(\frac{b}{\epsilon}\right )^q}{q},
\end{array}
$$
where $\epsilon>0$, $p>1$, $q>1$ and $\frac{1}{p}+\frac{1}{q}=1$.
\end{lemma}

\begin{lemma}(Gronwall's Inequality \cite{QSS}) \label{Lem.3.2}  Assume that $\{k_n\}$ and  $\{p_n\}$ are nonnegative sequences, and
the sequence $\{\phi_n\}$ satisfies
$$
\begin{array}{lll}
\displaystyle
\phi_0\leq q_0,\;\;\phi_n\leq q_0+\sum_{\ell=0}^{n-1}p_{\ell}+\sum_{\ell=0}^{n-1}k_{\ell}\phi_{\ell},\;n\geq1,
\end{array}
$$
where $q_{0}\geq0$. Then the sequence $\{\phi_n\}$ satisfies
$$
\begin{array}{lll}
\displaystyle
\phi_n\leq \left(q_0+\sum_{\ell=0}^{n-1}p_{\ell}\right)\exp\left(\sum_{\ell=0}^{n-1}k_{\ell}\right),\;n\geq1.
\end{array}
$$
\end{lemma}

\begin{lemma}\cite{GS} \label{Lem.3.3}  For any grid function $\bm{u}\in\mathcal{S}_h$, one has
$$
\begin{array}{lll}
\displaystyle
\frac{2}{3}||\bm{u}||^2\leq(\mathscr{L}\bm{u},\bm{u})\leq||\bm{u}||^2.
\end{array}
$$
\end{lemma}

\begin{lemma} \label{Lem.3.4}For any grid function $\bm{u}\in\mathcal{S}_h$, there exists a symmetric
positive difference operator denoted by $\mathscr{L}^{\frac{1}{2}}$,  such that
$$
\begin{array}{lll}
\displaystyle
(\mathscr{L}\bm{u},\bm{u})=(\mathscr{L}^{\frac{1}{2}}\bm{u},\mathscr{L}^{\frac{1}{2}}\bm{u}).
\end{array}
$$
\end{lemma}

\begin{proof}
Obviously, the corresponding matrix of operator $\mathscr{L}$ is given by
$$\mathbf{D}=
\left(
  \begin{array}{ccccc}
   \frac{5}{6} & \frac{1}{12}&0& \cdots & 0 \vspace{0.2 cm}\\
    \frac{1}{12} & \frac{5}{6} & \frac{1}{12} & 0 & \cdots \vspace{0.2 cm}\\
    \vdots &\vdots & \ddots & \ddots & \ddots \vspace{0.2 cm}\\
  0 &0& \ldots &\frac{5}{6}& \frac{1}{12}\vspace{0.2 cm}\\
 0& 0& \ldots & \frac{1}{12}& \frac{5}{6} \vspace{0.2 cm}\\
  \end{array}
\right).
$$
Obviously, $\mathbf{D}$ is real symmetric and positive definite. Hence, there exists
 a symmetric positive matrix denoted by $\mathbf{D}^{\frac{1}{2}}$
 such that
$$
\begin{array}{lll}
\displaystyle
(\mathscr{L}\bm{u},\bm{u})=h\bm{u}^{T}\mathbf{D}\bm{u}=h\mathbf{u}^{T}\left(\mathbf{D}^{\frac{1}{2}}\right)^2\bm{u}
=h\left(\mathbf{D}^{\frac{1}{2}}\bm{u}\right)^{T}\left(\mathbf{D}^{\frac{1}{2}}\bm{u}\right)
=(\mathscr{L}^{\frac{1}{2}}\bm{u},\mathscr{L}^{\frac{1}{2}}\bm{u}),
\end{array}
$$
where $\mathscr{L}^{\frac{1}{2}}$ is the associate operator of matrix $\mathbf{D}^{\frac{1}{2}}$.
Therefore, the proof is ended.
\end{proof}

\begin{lemma} \label{Lem.3.5} For any mesh function $\{\bm{u}^k|k=0,1,\ldots,n-1\}\in\mathcal{S}_h$, it holds that
$$
\begin{array}{lll}
\displaystyle
\sum_{k=0}^{n-1}\left(\delta_{0,t}^{1-\alpha}\bm{u}^{k},\bm{u}^{k+1}+\bm{u}^k\right)\geq0
\end{array}
$$
\end{lemma}
for $0<\alpha<1$, where $\alpha$ denotes $\alpha_1$ or $\alpha_2$.

\begin{proof}
Firstly, we have
$$
\begin{array}{lll}
\displaystyle
\sum_{k=0}^{n-1}\left(\delta_{0,t}^{1-\alpha}\bm{u}^{k},\bm{u}^{k+1}+\bm{u}^k\right)=h
\sum_{j=1}^{M-1}\sum_{k=0}^{n-1}\left(\sum_{\ell=0}^{k}\varpi_{\ell}^{(1-\alpha)}u_j^{k-\ell}\right)
\left(u_j^{k+1}+u_j^k\right)\vspace{0.2 cm}\\
\displaystyle
=h\sum_{j=1}^{M-1}\sum_{k=0}^{n-1}\left(\sum_{\ell=0}^{k}\varpi_{\ell}^{(1-\alpha)}u_j^{k-\ell}\right)u_j^k+
h\sum_{j=1}^{M-1}\sum_{k=0}^{n-1}\left(\sum_{\ell=0}^{k}\varpi_{\ell-1}^{(1-\alpha)}u_j^{k+1-\ell}
\right)u_j^{k+1}\vspace{0.2 cm}\\
\displaystyle
+h\sum_{j=1}^{M-1}\sum_{k=0}^{n-1}\left(
\varpi_{k}^{(1-\alpha)}u_j^0-\varpi_{-1}^{(1-\alpha)}u_j^{k+1}\right)u_j^{k+1}.
\end{array}
$$
 Letting $\varpi_{-1}^{(1-\alpha)}=0$ and noticing $u_j^0=0$, one has
$$
\begin{array}{lll}
\displaystyle
\sum_{k=0}^{n-1}\left(\delta_{0,t}^{1-\alpha}\bm{u}^{k},\mathbf{u}^{k+1}+\bm{u}^k\right)=
\vspace{0.2 cm}\\
\displaystyle
=h\sum_{j=1}^{M-1}\sum_{k=0}^{n-1}\left(\sum_{\ell=0}^{k}\varpi_{\ell}^{(1-\alpha)}u_j^{k-\ell}\right)u_j^k+
h\sum_{j=1}^{M-1}\sum_{k=0}^{n-1}\left(\sum_{\ell=0}^{k}\varpi_{\ell-1}^{(1-\alpha)}u_j^{k+1-\ell}
\right)u_j^{k+1}\vspace{0.2 cm}\\
=(\mathbf{A}_\alpha\bm{u},\bm{u})+(\mathbf{B}_\alpha\bm{u},\bm{u})=(\mathbf{E}_\alpha\bm{u},\bm{u}),
\end{array}
$$
where matrices $\mathbf{A}_\alpha$ and $\mathbf{B}_\alpha$ are
$$\mathbf{A}_\alpha=
\left(
  \begin{array}{ccccc}
   \varpi_{0}^{(1-\alpha)} & \frac{1}{2}\varpi_{1}^{(1-\alpha)} & \frac{1}{2}\varpi_{2}^{(1-\alpha)} & \cdots & \frac{1}{2}\varpi_{n}^{(1-\alpha)} \vspace{0.2 cm}\\
    \frac{1}{2}\varpi_{1}^{(1-\alpha)} & \varpi_{0}^{(1-\alpha)} & \frac{1}{2}\varpi_{1}^{(1-\alpha)} & \frac{1}{2}\varpi_{2}^{(1-\alpha)}& \cdots \vspace{0.2 cm}\\
    \vdots &\vdots & \ddots & \ddots & \ddots \vspace{0.2 cm}\\
  \frac{1}{2}\varpi_{n-1}^{(1-\alpha)} &\frac{1}{2}\varpi_{n-2}^{(1-\alpha)}& \ldots &\varpi_{0}^{(1-\alpha)}& \frac{1}{2}\varpi_{1}^{(1-\alpha)} \vspace{0.2 cm}\\
  \frac{1}{2}\varpi_{n}^{(1-\alpha)}& \frac{1}{2}\varpi_{n-1}^{(1-\alpha)}& \ldots & \frac{1}{2}\varpi_{1}^{(1-\alpha)}& \varpi_{0}^{(1-\alpha)} \vspace{0.2 cm}\\
  \end{array}
\right),
$$
and
$$\mathbf{B}_\alpha=
\left(
  \begin{array}{ccccc}
   \varpi_{-1}^{(1-\alpha)} & \frac{1}{2}\varpi_{0}^{(1-\alpha)} & \frac{1}{2}\varpi_{1}^{(1-\alpha)} & \cdots & \frac{1}{2}\varpi_{n-1}^{(1-\alpha)} \vspace{0.2 cm}\\
    \frac{1}{2}\varpi_{0}^{(1-\alpha)} & \varpi_{-1}^{(1-\alpha)} & \frac{1}{2}\varpi_{0}^{(1-\alpha)} & \frac{1}{2}\varpi_{1}^{(1-\alpha)}& \cdots \vspace{0.2 cm}\\
    \vdots &\vdots & \ddots & \ddots & \ddots \vspace{0.2 cm}\\
  \frac{1}{2}\varpi_{n-2}^{(1-\alpha)} &\frac{1}{2}\varpi_{n-3}^{(1-\alpha)}& \ldots &\varpi_{-1}^{(1-\alpha)}& \frac{1}{2}\varpi_{0}^{(1-\alpha)} \vspace{0.2 cm}\\
  \frac{1}{2}\varpi_{n-1}^{(1-\alpha)}& \frac{1}{2}\varpi_{n-2}^{(1-\alpha)}& \ldots & \frac{1}{2}\varpi_{0}^{(1-\alpha)}& \varpi_{-1}^{(1-\alpha)} \vspace{0.2 cm}\\
  \end{array}
\right).
$$

By the Grenander-Szeg\"{o} Theorem \cite{J}, if the generating function of matrix $\mathbf{E}_\alpha$
is nonnegative, then matrix $\mathbf{E}_\alpha$ is positive semi-definite. So, we only consider
 the generating function of matrix $\mathbf{E}_\alpha$
which is
$$
\begin{array}{lll}
\displaystyle
G_{\mathbf{E}_\alpha}(x,\alpha)&=&\displaystyle G_{\mathbf{A}_\alpha}(x,\alpha)+G_{\mathbf{B}_\alpha}(x,\alpha)\vspace{0.2 cm}\\
&=&\displaystyle
\varpi_{0}^{(1-\alpha)}+\frac{1}{2}\sum_{k=1}^{\infty}\varpi_{k}^{(1-\alpha)}\textmd{e}^{\textbf{i}kx}
+\frac{1}{2}\sum_{k=1}^{\infty}\varpi_{k}^{(1-\alpha)}\textmd{e}^{-\textbf{i}kx}
\vspace{0.2 cm}\\
&&\displaystyle
+\varpi_{-1}^{(1-\alpha)}+\frac{1}{2}\sum_{k=1}^{\infty}\varpi_{k-1}^{(1-\alpha)}\textmd{e}^{\textbf{i}kx}
+\frac{1}{2}\sum_{k=1}^{\infty}\varpi_{k-1}^{(1-\alpha)}\textmd{e}^{-\textbf{i}kx}
\vspace{0.2 cm}\\
&=&\displaystyle
\frac{1}{2}\sum_{k=0}^{\infty}\varpi_{k}^{(1-\alpha)}\textmd{e}^{\textbf{i}kx}
+\frac{1}{2}\sum_{k=0}^{\infty}\varpi_{k}^{(1-\alpha)}\textmd{e}^{-\textbf{i}kx}
+\frac{1}{2}\sum_{k=0}^{\infty}\varpi_{k}^{(1-\alpha)}\textmd{e}^{\textbf{i}(k+1)x}\vspace{0.2 cm}\\
&&\displaystyle+\frac{1}{2}\sum_{k=0}^{\infty}\varpi_{k}^{(1-\alpha)}\textmd{e}^{-\textbf{i}(k+1)x}
\vspace{0.2 cm}\\
&=&\displaystyle 2\left(\frac{4-3\alpha}{2(1-\alpha)}\right)^{1-\alpha}\left(2\sin\frac{x}{2}\right)^{1-\alpha}
\left(a^2+b^2\right)^{\frac{1-\alpha}{2}}\cos\left(\frac{x}{2}\right)
Z(x,\alpha),
\end{array}
$$
where
$$
\begin{array}{lll}
\displaystyle
a=1-\frac{2-\alpha}{4-3\alpha}\cos(x),\;\;
b=\frac{2-\alpha}{4-3\alpha}\sin(x),\vspace{0.2 cm}\\\displaystyle
Z(x,\alpha)=\cos\left(\frac{x}{2}+(1-\alpha)\left(\theta+\frac{x-\pi}{2}\right)\right),\;\;\;\theta=-\arctan\frac{b}{a}.
\end{array}
$$

Because $G_{\mathbf{E}_\alpha}(x,\alpha)$ is a real value and even function, we only consider
 the case of $x\in[0, \pi]$ for $G_{\mathbf{E}_\alpha}(x,\alpha)$. From the above formula, it is easy to know that only function
$Z(x,\alpha)$ needs to be studied for $0<\alpha<1$.

 Let
 $$
\begin{array}{lll}
\displaystyle
q(x,\alpha)=\frac{x}{2}+(1-\alpha)\left(\theta+\frac{x-\pi}{2}\right),\;x\in[0,\pi],\;\alpha\in(0,1).
\end{array}
$$
 One has
  $$
\begin{array}{lll}
\displaystyle
\frac{\partial q(x,\alpha)}{\partial x}
=\frac{4(2-\alpha)(1-\alpha)^2+2(4-3\alpha)(2-\alpha)^2\sin^2\left(\frac{x}{2}\right)}
{\left[(4-3\alpha)-(2-\alpha)\cos(x)\right]^2+\left[(2-\alpha)\sin(x)\right]^2}\geq0
\end{array}
$$
 which implies that $q(x,\alpha)$ is an increasing function with respect to $x$. And
 $$
\begin{array}{lll}
\displaystyle
q_{\min}(x,\alpha)=q(0,\alpha)=-\frac{\pi}{2}(1-\alpha),\;\;q_{\max}(x,\alpha)=q(\pi,\alpha)=\frac{\pi}{2}.
\end{array}
$$
 Hence,
  $$
\begin{array}{lll}
\displaystyle
Z(x,\alpha)=\cos\left(q(x,\alpha)\right)\geq0,
\end{array}
$$
 it immediately follows that
   $$
\begin{array}{lll}
\displaystyle
G_{\mathbf{E}_\alpha}(x,\alpha)\geq0.
\end{array}
$$
So the proof is completed.
\end{proof}

In the following, the first step is to prove the solvability of finite difference scheme (11), together with (12) and (13).

\begin{theorem}
The finite difference scheme (11), together with (12) and (13) is uniquely solvable.
\end{theorem}

\begin{proof}
Consider the homogeneous form of system (11),
$$
\begin{array}{lll}
\displaystyle
\mathscr{L}u_j^{k+1}=0.
\end{array}\eqno(14)
$$
Taking the inner product of (11) with $u_j^{k+1}$ gives
$$
\begin{array}{lll}
\displaystyle
\left(\mathscr{L}\bm{u}^{k+1},\bm{u}^{k+1}\right)=0.
\end{array}
$$

It follows from Lemma \ref{Lem.3.3} that
$$
\begin{array}{lll}
\displaystyle
||\bm{u}^{k+1}||^2=0,
\end{array}
$$
that is, ${u}_j^{k+1}=0$. This finishes the proof.
\end{proof}

Now, we give the stability result.

\begin{theorem} The finite difference scheme (11), together with (12) and (13) is unconditionally
 stable with respect to the initial value.
\end{theorem}

\begin{proof}
Let $u_j^k$ and $v_j^k$ be the solutions of the following two equations, respectively,
$$
\begin{array}{lll}
\displaystyle
\mathscr{L}\left(u_j^{k+1}-u_j^{k}\right)=\tau K_1\delta_{0,t}^{1-\alpha_1}\delta_x^2u_j^{k}
-\tau K_2\delta_{0,t}^{1-\alpha_2}\mathscr{L}u_j^{k}  +\tau\mathscr{L}f_j^{k+\frac{1}{2}},\vspace{0.2 cm}\\
\hspace{5.2 cm}\;\;1\leq j\leq M-1,\;1\leq k\leq N,\vspace{0.2 cm}\\
\displaystyle
u_j^0=0,\;\;0\leq j\leq M,
\vspace{0.2 cm}\\
\displaystyle
u_0^k=\varphi_1(t_k),\;u_M^k=\varphi_2(t_k),\;\;1\leq k\leq N;
\end{array}
$$
and
$$
\begin{array}{lll}
\displaystyle
\mathscr{L}\left(v_j^{k+1}-v_j^{k}\right)=\tau K_1\delta_{0,t}^{1-\alpha_1}\delta_x^2v_j^{k}
-\tau K_2\delta_{0,t}^{1-\alpha_2}\mathscr{L}v_j^{k}  +\tau\mathscr{L}f_j^{k+\frac{1}{2}},\vspace{0.2 cm}\\
\hspace{5.2 cm}\;\;1\leq j\leq M-1,\;1\leq k\leq N,\vspace{0.2 cm}\\
\displaystyle
v_j^0=\rho_j,\;\;0\leq j\leq M,
\vspace{0.2 cm}\\
\displaystyle
v_0^k=\varphi_1(t_k),\;v_M^k=\varphi_2(t_k),\;\;1\leq k\leq N.
\end{array}
$$
Denote $\varepsilon_j^k=v_j^k-u_j^k$. Then one can gets
the following perturbation equation,
$$
\begin{array}{lll}
\displaystyle
\mathscr{L}\left(\varepsilon_j^{k+1}-\varepsilon_j^{k}\right)=\tau K_1\delta_{0,t}^{1-\alpha_1}\delta_x^2\varepsilon_j^{k}
-\tau K_2\delta_{0,t}^{1-\alpha_2}\mathscr{L}\varepsilon_j^{k}, 1\leq j\leq M-1,\;1\leq k\leq N,\vspace{0.2 cm}\\
\displaystyle
\varepsilon_j^0=\rho_j,\;\;0\leq j\leq M,
\vspace{0.2 cm}\\
\displaystyle
\varepsilon_0^k=0,\;\varepsilon_M^k=0,\;\;1\leq k\leq N.
\end{array}\eqno(15)
$$

Taking the discrete inner product with $(\bm{\varepsilon}^{k+1}+\bm{\varepsilon}^k)$
on the both sides of the equation (15)
and summing up for $k$ from $0$ to $n-1$ give
$$
\begin{array}{lll}
\displaystyle
\sum_{k=0}^{n-1}\left(\mathscr{L}\left(\bm{\varepsilon}^{k+1}-\bm{\varepsilon}^{k}\right),(\bm{\varepsilon}^{k+1}+\bm{\varepsilon}^k)\right)=&\displaystyle\tau K_1\sum_{k=0}^{n-1}\left(\delta_{0,t}^{1-\alpha_1}\delta_x^2\bm{\varepsilon}^{k},(\bm{\varepsilon}^{k+1}+\bm{\varepsilon}^k)\right)\vspace{0.2 cm}\\&\displaystyle
-\tau K_2\sum_{k=0}^{n-1}\left(\delta_{0,t}^{1-\alpha_2}\mathscr{L}\bm{\varepsilon}^{k},(\bm{\varepsilon}^{k+1}+\bm{\varepsilon}^k)\right).
\end{array}
$$
Furthermore, one has the following results
$$
\begin{array}{lll}
\displaystyle
\sum_{k=0}^{n-1}\left(\mathscr{L}\left(\bm{\varepsilon}^{k+1}-\bm{\varepsilon}^{k}\right),(\bm{\varepsilon}^{k+1}+\bm{\varepsilon}^k)\right)=
\sum_{k=0}^{n-1}\left[\left(\mathscr{L}\bm{\varepsilon}^{k+1},\bm{\varepsilon}^{k+1}\right)
-\left(\mathscr{L}\bm{\varepsilon}^{k},\bm{\varepsilon}^{k}\right)\right],
\end{array}
$$
$$
\begin{array}{lll}
\displaystyle
\sum_{k=0}^{n-1}\left(\delta_{0,t}^{1-\alpha_1}\delta_x^2\bm{\varepsilon}^{k},(\bm{\varepsilon}^{k+1}+\bm{\varepsilon}^k)\right)
=-\sum_{k=0}^{n-1}\left(\delta_{0,t}^{1-\alpha_1}\delta_x\bm{\varepsilon}^{k},\delta_x(\bm{\varepsilon}^{k+1}+\bm{\varepsilon}^k)\right)\leq0,
\end{array}
$$
and
$$
\begin{array}{lll}
\displaystyle
-\sum_{k=0}^{n-1}\left(\delta_{0,t}^{1-\alpha_2}\mathscr{L}\bm{\varepsilon}^{k},(\bm{\varepsilon}^{k+1}+\bm{\varepsilon}^k)\right)=-
\sum_{k=0}^{n-1}\left(\delta_{0,t}^{1-\alpha_2}\mathscr{L}^{\frac{1}{2}}\bm{\varepsilon}^{k},
\mathscr{L}^{\frac{1}{2}}(\bm{\varepsilon}^{k+1}+\bm{\varepsilon}^k)\right)\leq0
\end{array}
$$
by follows that Lemmas \ref{Lem.3.4} and \ref{Lem.3.5}. Hence,
$$
\begin{array}{lll}
\displaystyle
\sum_{k=0}^{n-1}\left[\left(\mathscr{L}\bm{\varepsilon}^{k+1},\bm{\varepsilon}^{k+1}\right)
-\left(\mathscr{L}\bm{\varepsilon}^{k},\bm{\varepsilon}^{k}\right)\right]\leq0,
\end{array}
$$
i.e.,
$$
\begin{array}{lll}
\displaystyle
\left(\mathscr{L}\bm{\varepsilon}^{n},\bm{\varepsilon}^{n}\right)
\leq\left(\mathscr{L}\bm{\varepsilon}^{0},\bm{\varepsilon}^{0}\right).
\end{array}
$$
Using Lemma \ref{Lem.3.3} again, one has,
$$
\begin{array}{lll}
\displaystyle
||\bm{\varepsilon}^{n}||\leq\frac{\sqrt{6}}{2}||\bm{\varepsilon}^{0}||=\frac{\sqrt{6}}{2}||\bm{\rho}||.
\end{array}
$$
This ends the proof.
\end{proof}

Finally, we study the convergence of the above compact difference scheme.

\begin{theorem}
Assume that the solution $u(x, t)$ of equation (6), together with (7) and (8) is suitably
smooth, and let $\{u_j^k|
0\leq j\leq M, 0 \leq k \leq N\}$ be the solution of the finite difference scheme
(11)--(13). Set $e_j^k=u(x_j, t_k)-u_j^k$, then
$$
\begin{array}{lll}
\displaystyle
||\bm{e}^{n}||\leq \frac{3}{2}C_1\exp(T)\sqrt{2LT}(\tau^2+h^4).
\end{array}
$$
\end{theorem}

\begin{proof}
From the above analysis, we have the following error system,
$$
\begin{array}{lll}
\displaystyle
\mathscr{L}\left(e_j^{k+1}-e_j^{k}\right)=\tau K_1\delta_{0,t}^{1-\alpha_1}\delta_x^2e_j^{k}
-\tau K_2\delta_{0,t}^{1-\alpha_2}\mathscr{L}e_j^{k}+\tau R_j^k,\vspace{0.2 cm}\\\;\;1\leq j\leq M-1,\;1\leq k\leq N,\vspace{0.2 cm}\\
\displaystyle
e_j^0=0,\;\;0\leq j\leq M,
\vspace{0.2 cm}\\
\displaystyle
e_0^k=0,\;e_M^k=0,\;\;1\leq k\leq N.
\end{array}\eqno(16)
$$
Taking the discrete inner product with $(\bm{e}^{k+1}+\bm{e}^k)$
on the both sides of equation (16)
and summing up for $k$ from $0$ to $n-1$, one gets
$$
\begin{array}{lll}
\displaystyle
\sum_{k=0}^{n-1}\left(\mathscr{L}\left(\bm{e}^{k+1}-\bm{e}^{k}\right),(\bm{e}^{k+1}+\bm{e}^k)\right)=\tau K_1\sum_{k=0}^{n-1}\left(\delta_{0,t}^{1-\alpha_1}
\delta_x^2\bm{e}^{k},(\bm{e}^{k+1}+\bm{\varepsilon}^k)\right)\vspace{0.2 cm}\\\displaystyle
-\tau K_2\sum_{k=0}^{n-1}\left(\delta_{0,t}^{1-\alpha_2}\mathscr{L}\bm{e}^{k},(\bm{e}^{k+1}+\bm{e}^k)\right)+
\tau\sum_{k=0}^{n-1}\left(\bm{R}^k,(\bm{e}^{k+1}+\bm{e}^k)\right).
\end{array}
$$
It follows from Lemma \ref{Lem.3.1} that
$$
\begin{array}{lll}
\displaystyle
\tau\sum_{k=0}^{n-1}\left(\bm{R}^k,(\bm{e}^{k+1}+\bm{e}^k)\right)\leq\tau\sum_{k=0}^{n-1}
||\bm{R}^k||\left(||\bm{e}^{k+1}||+||\bm{e}^{k}||\right)\vspace{0.2 cm}\\ \displaystyle \leq
\frac{3}{2}\tau\sum_{k=0}^{n-1}
||\bm{R}^k||^2+\frac{1}{3}\tau\sum_{k=0}^{n-1}
\left(||\bm{e}^{k+1}||^2+||\bm{e}^{k}||^2\right).
\end{array}
$$
Hence, one further arrives at
$$
\begin{array}{lll}
\displaystyle
\frac{2}{3}||\bm{e}^{n}||^2\leq
\frac{3}{2}\tau\sum_{k=0}^{n-1}
||\bm{R}^k||^2+\frac{1}{3}\tau\sum_{k=0}^{n-1}
\left(||\bm{e}^{k+1}||^2+||\bm{e}^{k}||^2\right)\vspace{0.2 cm}\\ \displaystyle \leq
\frac{1}{3}||\bm{e}^{n}||^2+\frac{2}{3}\tau\sum_{k=0}^{n-1}
||\bm{e}^{k}||^2+\frac{3}{2}LTC_1^2(\tau^2+h^4)^2,
\end{array}
$$
i.e.,
$$
\begin{array}{lll}
\displaystyle
||\bm{e}^{n}||^2\leq
2\tau\sum_{k=0}^{n-1}
||\bm{e}^{k}||^2+\frac{9}{2}LTC_1^2(\tau^2+h^4)^2,
\end{array}
$$
Applying Lemma \ref{Lem.3.2} to it yields
$$
\begin{array}{lll}
\displaystyle
||\bm{e}^{n}||\leq \frac{3}{2}C_1\exp(T)\sqrt{2LT}(\tau^2+h^4).
\end{array}
$$
The proof is thus finished.
\end{proof}

\section{Application to the fractional Cable equation in two space dimensions}

 In this section, we study the following two-dimensional Cable equation
$$
\begin{array}{ll}
\displaystyle
\frac{\partial u(x,y,t)}{\partial t}=&\displaystyle K_1\,_{\mathrm{RL}}{{{\mathrm{D}}}}_{0,t}^{1-\alpha_1}\left(
\frac{\partial^2u(x,y,t)}{\partial x^2}+\frac{\partial^2u(x,y,t) }{\partial y^2}
\right)\vspace{0.2 cm}\\&\displaystyle
- K_2\,_{\mathrm{RL}}{{{\mathrm{D}}}}_{0,t}^{1-\alpha_2}\left( u(x,y,t)\right) +f(x,y,t),
\vspace{0.2 cm}\\&\hspace{2.5 cm}\;(x,y)\in\Omega=(0,L_x)\times(0,L_y),\;t\in(0,T],
\end{array}\eqno(17)
$$
with initial condition
$$
\begin{array}{ll}
\displaystyle
u(x,y,0)=0,\;\;(x,y)\in\overline{\Omega}=[0,L_x]\times[0,L_y],
\end{array}\eqno(18)
$$
and boundary value conditions
$$
\begin{array}{ll}
\displaystyle
u(x,y,t)=\varphi(x,y,t),\;\;  (x,y)\in{\partial\Omega},
\end{array}\eqno(19)
$$
where $0<\alpha_1,\alpha_2<1$, $K_1$ and $K_2$ are two constants,
$f(x,y, t)$ and $\varphi(x,y,t)$ are given suitably
smooth functions.

\subsection{Development of numerical algorithm}

Denote $x_i=i h_x,\;0\leq i\leq M_1$, $y_j=j h_y,\;0\leq j\leq M_2$, where $h_x=L_x/M_1$, $h_y=L_y/M_2$ and $M_1,M_2$
are two positive integers. Define $\overline{\Omega}_{h_{x}h_{y}}=\{(x_i,y_j)|0\leq i\leq M_1,0\leq j\leq M_2\}$,
${\Omega}_{h_{x}h_{y}}=\overline{\Omega}_{h_{x}h_{y}}\cap \Omega$ and $\partial{\Omega}_{h_{x}h_{y}}=\overline{\Omega}_{h_{x}h_{y}}\cap \partial\Omega$.
For any grid function $\bm{u}\in V_{h_{x}h_{y}}=\{\bm{u}|\bm{u}=\{u_{i,j}\}, 0\leq i\leq M_1, 0\leq j\leq M_2\}$, define the difference operators as,
$$
\begin{array}{ll}
\displaystyle
\delta_x u_{i-\frac{1}{2},j}=\frac{1}{h_x}\left(u_{i,j}-u_{i-1,j}\right),\;\;
\delta_x^2 u_{i,j}=\frac{1}{h_x}\left(\delta_x u_{i+\frac{1}{2},j}
-\delta_x u_{i-\frac{1}{2},j}\right),\vspace{0.2 cm}\\ \displaystyle
\delta_y u_{i,j-\frac{1}{2}}=\frac{1}{h_y}\left(u_{i,j}-u_{i,j-1}\right),
\;\;\delta_y^2 u_{i,j}=\frac{1}{h_y}\left(\delta_y u_{i,j+\frac{1}{2}}-\delta_y u_{i,j-\frac{1}{2}}\right),
\end{array}
$$
and the spatial compact difference operators as,
$$\mathscr{L}_xu_{i,j}=\left\{
\begin{array}{lll}
\left(I+\frac{h_x^2}{12}\delta_x^2\right) u_{i,j},\;\;1\leq i\leq M_1-1,\;0\leq j\leq M_2,\vspace{0.2 cm}\\
u_{i,j},\;\;i=0\; \textrm{or}\; M_1,\;0\leq j\leq M_2,
\end{array}\right.
$$
$$\mathscr{L}_yu_{i,j}=\left\{
\begin{array}{lll}
\left(I+\frac{h_y^2}{12}\delta_y^2\right) u_{i,j},\;\;1\leq j\leq M_2-1,\;0\leq i\leq M_1,\vspace{0.2 cm}\\
u_{i,j},\;\;j=0\; \textrm{or}\; M_2,\;0\leq i\leq M_1.
\end{array}\right.
$$

Now we consider equation (17) at point $\left(x_i,y_j,t_{k+\frac{1}{2}}\right)$. Operating the operator $\mathscr{L}_x\mathscr{L}_y$
on both sides of it leads to,
$$
\begin{array}{ll}
\displaystyle
\mathscr{L}_x\mathscr{L}_y\frac{\partial u\left(x_i,y_j,t_{k+\frac{1}{2}}\right)}{\partial t}=K_1\,_{\mathrm{RL}}{{{\mathrm{D}}}}_{0,t}^{1-\alpha_1}\left(\mathscr{L}_y
\mathscr{L}_x\frac{\partial^2u\left(x_i,y_j,t_{k+\frac{1}{2}}\right) }{\partial x^2}\right.\vspace{0.2 cm}\\\left.\displaystyle
+\mathscr{L}_x\mathscr{L}_y\frac{\partial^2u\left(x_i,y_j,t_{k+\frac{1}{2}}\right) }{\partial y^2}
\right)
- K_2\,_{\mathrm{RL}}{{{\mathrm{D}}}}_{0,t}^{1-\alpha_2}\mathscr{L}_x\mathscr{L}_y u\left(x_i,y_j,t_{k+\frac{1}{2}}\right)
\vspace{0.2 cm}\\\displaystyle +\mathscr{L}_x\mathscr{L}_yf\left(x_i,y_j,t_{k+\frac{1}{2}}\right),
1\leq i\leq M_1-1,\;1\leq j\leq M_2-1,\;0\leq k\leq N-1.
\end{array}
$$
Similar to the one-dimensional case, one easily has
$$
\begin{array}{ll}
\displaystyle
\mathscr{L}_x\mathscr{L}_y
\left(u(x_i,y_j,t_{k+1})-u(x_i,y_j,t_{k})\right)\vspace{0.2 cm}\\\displaystyle
=\tau K_1\delta_{0,t}^{1-\alpha_1}\left(\mathscr{L}_y
\delta_x^2+\mathscr{L}_x\delta_y^2
\right)u\left(x_i,y_j,t_{k}\right)
- \tau K_2\delta_{0,t}^{1-\alpha_2}\mathscr{L}_x\mathscr{L}_y u\left(x_i,y_j,t_{k}\right)
\vspace{0.2 cm}\\\displaystyle +\tau\mathscr{L}_x\mathscr{L}_yf\left(x_i,y_j,t_{k+\frac{1}{2}}\right)+\tau R_{i,j}^k,
\vspace{0.2 cm}\\\hspace{3.0 cm}\;1\leq i\leq M_1-1,\;1\leq j\leq M_2-1,\;0\leq k\leq N-1,
\end{array}\eqno(20)
$$
where there exists a positive constant $C_2$ such that $|R_{i,j}^k|\leq C_2(\tau^2+h_x^4+h_y^4)$.

Omitting the small term $R_{i,j}^k$ in (20) and replacing the grid function $u\left(x_i,y_j,t_{k}\right)$ with its
numerical approximation $u_{i,j}^k$, one gets a high-order compact scheme in the following form,
$$
\begin{array}{lll}
\displaystyle
\mathscr{L}_x\mathscr{L}_y
\left(u_{i,j}^{k+1}-u_{i,j}^k\right)
=\tau K_1\delta_{0,t}^{1-\alpha_1}\left(\mathscr{L}_y
\delta_x^2+\mathscr{L}_x\delta_y^2
\right)u_{i,j}^k-\vspace{0.2 cm}\\ \tau K_2\delta_{0,t}^{1-\alpha_2}\mathscr{L}_x\mathscr{L}_y u_{i,j}^k
\displaystyle
 +\tau\mathscr{L}_x\mathscr{L}_yf_{i,j}^{k+\frac{1}{2}}+\tau R_{i,j}^k,\vspace{0.2 cm}\\
\;1\leq i\leq M_1-1,\;1\leq j\leq M_2-1,\;0\leq k\leq N-1,\vspace{0.2 cm}\\
\end{array}\eqno(21)
$$
$$
\begin{array}{lll}
\displaystyle
\displaystyle
u_{i,j}^0=0,\;\;\;\;\;0\leq i\leq M_1,\;0\leq j\leq M_2,\vspace{0.2 cm}\\
\end{array}\eqno(22)
$$
$$
\begin{array}{lll}
\displaystyle
\displaystyle
u_{i,j}^k=\varphi\left(x_i,y_j,t_{k}\right),\;\;\;\;(x_i,y_j)\in\partial{\Omega}_{h_{x}h_{y}}.
\end{array}\eqno(23)
$$

\subsection{Solvability, stability and convergence analysis}
Define
$$
\begin{array}{ll}
\displaystyle
\mathcal{S}_{h_xh_y}=\{\bm{u}|\bm{u}=\{u_{i,j}\},u|_{\partial{\Omega}_{h_{x}h_{y}}}=0\}.
\end{array}
$$
Then for any $\bm{u},\bm{v}\in \mathcal{S}_{h_xh_y}$,
define the inner products as
$$
\begin{array}{ll}
\displaystyle
(\bm{u},\bm{v})=h_xh_y\sum_{i=1}^{M_1-1}\sum_{j=1}^{M_2-1}u_{i,j}v_{i,j},\;\;\;\vspace{0.2 cm}\\\displaystyle
(\delta_x\delta_y\bm{u},\delta_x\delta_y\bm{v})=h_xh_y\sum_{i=1}^{M_1}\sum_{j=1}^{M_2}
\left(\delta_x\delta_yu_{i-{\frac{1}{2}},j-{\frac{1}{2}}}\right)\left(\delta_x\delta_yv_{i-{\frac{1}{2}},j-{\frac{1}{2}}}\right),\vspace{0.2 cm}\\\displaystyle
(\delta_x \bm{u},\delta_x \bm{v})=h_xh_y\sum_{i=1}^{M_1}\sum_{j=1}^{M_2-1}
\left(\delta_x u_{i-{\frac{1}{2}},j }\right)\left(\delta_x v_{i-{\frac{1}{2}},j }\right),\vspace{0.2 cm}\\\displaystyle
(\delta_y \bm{u},\delta_y \bm{v})=h_xh_y\sum_{i=1}^{M_1-1}\sum_{j=1}^{M_2}
\left(\delta_y u_{i,j-{\frac{1}{2}} }\right)\left(\delta_y v_{i,j-{\frac{1}{2}},j }\right).
\end{array}
$$

We temporarily leave to give several definition(s) and lemmas \cite{L} which will be utilized later on.
\begin{definition}
If $\mathbf{A}=(a_{ij})$ is an $m \times n$ matrix and
$\mathbf{B}=(b_{ij})$ is a $p \times q$ matrix, then the Kronecker product $\mathbf{A}\otimes \mathbf{B}$ is an $mp \times nq$ block matrix
and denoted by
$$\mathbf{A}\otimes \mathbf{B}=
\left(
  \begin{array}{cccc}
   a_{11}\mathbf{B} & a_{12}\mathbf{B}&  \cdots & a_{1n}\mathbf{B} \vspace{0.2 cm}\\
     a_{21}\mathbf{B} &   a_{22}\mathbf{B} & \cdots&  a_{2n}\mathbf{B} \vspace{0.2 cm}\\
    \vdots &\vdots  & \ddots & \vdots \vspace{0.2 cm}\\
 a_{m1}\mathbf{B}& a_{m2}\mathbf{B} &  \cdots& a_{mn}\mathbf{B} \vspace{0.2 cm}\\
  \end{array}
\right).
$$
\end{definition}

\begin{lemma}\label{Lem.4.1}
Suppose that $\mathbf{A}\in \mathds{R}^{n\times n}$ has eigenvalues $\{\lambda_j\}^n_{j=1}$, and $\mathbf{B}\in \mathds{R}^{m\times m}$ has eigenvalues $\{\mu_j\}^m_{j=1}$.
Then the $mn$ eigenvalues of $\mathbf{A}\otimes \mathbf{B}$ are given below,
$$
\begin{array}{ll}
\displaystyle
 \lambda_1\mu_1,\ldots,\lambda_1\mu_m, \lambda_2\mu_1,\ldots, \lambda_2\mu_m,\ldots, \lambda_n\mu_1,\ldots,\lambda_n\mu_m.
\end{array}
$$
\end{lemma}

\begin{lemma}\label{Lem.4.2}
Let $\mathbf{A}\in \mathds{R}^{m\times n}$, $\mathbf{B}\in \mathds{R}^{r\times s}$,
$\mathbf{C}\in \mathds{R}^{n\times p}$, $\mathbf{D}\in \mathds{R}^{s\times t}$. Then
$$
\begin{array}{ll}
\displaystyle
 \left(\mathbf{A}\otimes \mathbf{B}\right)\left(\mathbf{C}\otimes \mathbf{D}\right)=\mathbf{AC}\otimes \mathbf{BD}.
\end{array}
$$
Moreover, if $\mathbf{A}\in \mathds{R}^{n\times n}$, $\mathbf{B}\in \mathds{R}^{m\times m}$, $\mathbf{I}_n$ and $\mathbf{I}_m$
are unit matrices of orders $n$, $m$, respectively, then matrices
 $\mathbf{I}_m\otimes \mathbf{A}$ and $\mathbf{B}\otimes\mathbf{I}_n$ can commute with each other.
\end{lemma}

\begin{lemma}\label{Lem.4.3}
For both $\mathbf{A}$ and $\mathbf{B}$, there holds
$$
\begin{array}{ll}
\displaystyle
\left(\mathbf{A}\otimes \mathbf{B}\right)^{T}=\mathbf{A}^{T}\otimes \mathbf{B}^{T}.
\end{array}
$$
\end{lemma}

Next we give an inequality for the mesh function.
\begin{lemma}\cite{C}\label{Lem.4.4}
For any mesh function $\bm{u}\in \mathcal{S}_{h_xh_y}$, there holds
$$
\begin{array}{ll}
\displaystyle
\frac{1}{3}||\bm{u}||^2\leq
\left(\mathscr{L}_x\mathscr{L}_y\bm{u},\bm{u}\right)\leq||\bm{u}||^2.
\end{array}
$$
\end{lemma}

Now we return to discuss the derived compact scheme.
\begin{lemma}\label{Lem.4.5}
For any mesh functions $\bm{u},\bm{v}\in \mathcal{S}_{h_xh_y}$, there has a symmetric positive
 definite operator $\mathscr{Q}_{x+y}$, such that
$$
\begin{array}{ll}
\displaystyle
((\mathscr{L}_y
\delta_x^2+\mathscr{L}_x\delta_y^2)\bm{u},\bm{v})=-(\mathscr{Q}_{x+y}\bm{u},\mathscr{Q}_{x+y}\bm{v}).
\end{array}
$$
\end{lemma}

\begin{proof}
Denote
$$
\begin{array}{ll}
\bm{u}=(u_{1,1},u_{2,1},\ldots,u_{M_1-1,1},\ldots,u_{1,M_2-1},u_{2,M_2-1},\ldots,u_{M_1-1,M_2-1})^{T}
\end{array}
$$
and
$$
\begin{array}{ll}
\bm{v}=(v_{1,1},v_{2,1},\ldots,v_{M_1-1,1},\ldots,v_{1,M_2-1},v_{2,M_2-1},\ldots,v_{M_1-1,M_2-1})^{T}.
\end{array}
$$
Then
$
((\mathscr{L}_y
\delta_x^2+\mathscr{L}_x\delta_y^2)\bm{u},\bm{v})
$
can be rewritten as
$$
\begin{array}{ll}
\displaystyle
((\mathscr{L}_y
\delta_x^2+\mathscr{L}_x\delta_y^2)\bm{u},\bm{v})=h_1h_2\bm{v}^{T}\mathbf{L}_{x+y}\bm{u},
\end{array}
$$
where
$$
\begin{array}{ll}
\displaystyle
\mathbf{L}_{x+y}=\frac{1}{h_1^2}\left(\mathbf{C}_2\otimes \mathbf{I}_1\right)\left(\mathbf{I}_2\otimes \mathbf{D}_1\right)
+\frac{1}{h_2^2}\left(\mathbf{I}_2\otimes \mathbf{C}_1\right)\left(\mathbf{D}_2\otimes \mathbf{I}_1\right).
\end{array}
$$
Here $\mathbf{I}_p$ are the unit matrices of order $M_p-1$ $(p=1,2)$,
$\mathbf{C}_p$ and $\mathbf{D}_p$ $(p=1,2)$ are defined by
$$\mathbf{C}_p=
\left(
  \begin{array}{ccccc}
   \frac{5}{6} & \frac{1}{12}&0& \cdots & 0 \vspace{0.2 cm}\\
    \frac{1}{12} & \frac{5}{6} & \frac{1}{12} & 0 & \cdots \vspace{0.2 cm}\\
    \vdots &\vdots & \ddots & \ddots & \ddots \vspace{0.2 cm}\\
  0 &0& \ldots &\frac{5}{6}& \frac{1}{12}\vspace{0.2 cm}\\
 0& 0& \ldots & \frac{1}{12}& \frac{5}{6} \vspace{0.2 cm}\\
  \end{array}
\right)_{(M_p-1)\times(M_p-1)}
$$
and
$$\mathbf{D}_p=
\left(
  \begin{array}{ccccc}
   -2& 1&0& \cdots & 0 \vspace{0.2 cm}\\
   1 & -2 & 1 & 0 & \cdots \vspace{0.2 cm}\\
    \vdots &\vdots & \ddots & \ddots & \ddots \vspace{0.2 cm}\\
  0 &0& \ldots &-2& 1\vspace{0.2 cm}\\
 0& 0& \ldots & 1& -2 \vspace{0.2 cm}\\
  \end{array}
\right)_{(M_p-1)\times(M_p-1)}.
$$

By using Lemmas \ref{Lem.4.2} and \ref{Lem.4.3} and the fact
$$
\begin{array}{lll}
\displaystyle
\mathbf{L}_{x+y}^{T}&=&\displaystyle\frac{1}{h_1^2}\left(\mathbf{I}_2\otimes \mathbf{D}_1\right)^{T}\left(\mathbf{C}_2\otimes \mathbf{I}_1\right)^{T}
+\frac{1}{h_2^2}\left(\mathbf{D}_2\otimes \mathbf{I}_1\right)^{T}\left(\mathbf{I}_2\otimes \mathbf{C}_1\right)^{T}\vspace{0.2 cm}\\
&=&\displaystyle\frac{1}{h_1^2}\left(\mathbf{I}_2\otimes \mathbf{D}_1\right)\left(\mathbf{C}_2\otimes \mathbf{I}_1\right)
+\frac{1}{h_2^2}\left(\mathbf{D}_2\otimes \mathbf{I}_1\right)\left(\mathbf{I}_2\otimes \mathbf{C}_1\right)\vspace{0.2 cm}\\
&=&\displaystyle\frac{1}{h_1^2}\left(\mathbf{C}_2\otimes \mathbf{I}_1\right)\left(\mathbf{I}_2\otimes \mathbf{D}_1\right)
+\frac{1}{h_2^2}\left(\mathbf{I}_2\otimes \mathbf{C}_1\right)\left(\mathbf{D}_2\otimes \mathbf{I}_1\right)
\vspace{0.2 cm}\\
&=&\displaystyle\mathbf{L}_{x+y},
\end{array}
$$
we know that matrix $\mathbf{L}_{x+y}$ is real symmetric.
It follows from Lemma \ref{Lem.4.1} that matrix $\mathbf{L}_{x+y}$ is negative definite since its
eigenvalues are all negative. Hence, we can declare that matrix $\mathbf{L}_{x+y}$
is real symmetric and negative definite, and there exist an orthogonal matrix $\mathbf{H}_{x+y}$ and a diagonal matrix $\mathbf{\Lambda}_{x+y}$ such that
$$
\begin{array}{lll}
\displaystyle
\mathbf{L}_{x+y}=-\mathbf{H}_{x+y}^{T}\mathbf{\Lambda}_{x+y}\mathbf{H}_{x+y}
=-\left(\mathbf{\Lambda}_{x+y}^{\frac{1}{2}}\mathbf{H}_{x+y}\right)^{T}\left(\mathbf{\Lambda}_{x+y}^{\frac{1}{2}}\mathbf{H}_{x+y}\right)
=-\mathbf{Q}_{x+y}^{T}\mathbf{Q}_{x+y}.
\end{array}
$$
Hence, we have
$$
\begin{array}{lll}
\displaystyle
((\mathscr{L}_y
\delta_x^2+\mathscr{L}_x\delta_y^2)\bm{u},\bm{v})&=&\displaystyle h_1h_2\bm{v}^{T}\mathbf{L}_{x+y}\bm{u}\vspace{0.2 cm}\\&=&\displaystyle
-h_1h_2\left(\mathbf{Q}_{x+y}\bm{v}\right)^{T}\mathbf{Q}_{x+y}\bm{u}=-(\mathscr{Q}_{x+y}\bm{u},\mathscr{Q}_{x+y}\bm{v}),
\end{array}
$$
where $\mathscr{Q}_{x+y}$ is a symmetric and positive definite operator,
 so is $\mathbf{Q}_{x+y}$. Therefore, the proof is ended.
\end{proof}

\begin{lemma}\cite{JS}\label{Lem.4.6}
For any mesh function $\bm{u}\in \mathcal{S}_{h_xh_y}$, there exists a symmetric positive
 definite operator $\mathscr{Q}_{xy}$ such that
$$
\begin{array}{ll}
\displaystyle
(\mathscr{L}_x\mathscr{L}_y
\bm{u},\bm{v})=(\mathscr{Q}_{xy}\bm{u},\mathscr{Q}_{xy}\bm{v}).
\end{array}
$$
\end{lemma}

In the following, we firstly give the solvability of finite difference scheme (21), together with (22) and (23).

\begin{theorem}
The finite difference scheme (21), together with (22) and (23) is uniquely solvable.
\end{theorem}
\begin{proof}
As the same as the one dimensional case, we easily get the homogeneous form of difference scheme (21) is
$$
\begin{array}{lll}
\displaystyle
\mathscr{L}_x\mathscr{L}_y
u_{i,j}^{k+1}
=0.
\end{array}
$$
Taking the inner product with $u_{i,j}^{k+1}$ and using Lemma \ref{Lem.4.4} lead to
$$
\begin{array}{ll}
\displaystyle
||\bm{u}^{k+1}||=0,
\end{array}
$$
i.e., ${u}^{k+1}_{i,j}=0$, which shows that the conclusion holds.
\end{proof}

Next, we study the stability result of finite difference scheme (21)--(23).

\begin{theorem}
The finite difference scheme (21), together with (22) and (23) is unconditionally stable with respect to the initial value.
\end{theorem}

\begin{proof} Due to the perturbation equation
$$
\begin{array}{ll}
\displaystyle
\mathscr{L}_x\mathscr{L}_y
\left(\varepsilon_{i,j}^{k+1}-\varepsilon_{i,j}^k\right)
=\tau K_1\delta_{0,t}^{1-\alpha_1}\left(\mathscr{L}_y
\delta_x^2+\mathscr{L}_x\delta_y^2
\right)\varepsilon_{i,j}^k- \tau K_2\delta_{0,t}^{1-\alpha_2}\mathscr{L}_x\mathscr{L}_y \varepsilon_{i,j}^k\vspace{0.2 cm}\\
\displaystyle
\;\hspace{3.2 cm}1\leq i\leq M_1-1,\;1\leq j\leq M_2-1,\;0\leq k\leq N-1,\vspace{0.2 cm}\\
\displaystyle
\varepsilon_{i,j}^0=\rho_{i,j},\;\;\;\;\;1\leq i\leq M_1,\;1\leq j\leq M_2,\vspace{0.2 cm}\\
\displaystyle
\varepsilon_{i,j}^k=0,\;\;\;\;(x_i,y_j)\in\partial{\Omega}_{h_{x}h_{y}},
\end{array}
$$
taking the inner product with $\left(\bm{\varepsilon}^{k+1}+\bm{\varepsilon}^k\right)$
for the first equation and summing up for $k$ from 0 to $n-1$ yield
$$
\begin{array}{ll}
\displaystyle
\sum_{k=0}^{n-1}\left(\mathscr{L}_x\mathscr{L}_y
\left(\bm{\varepsilon}^{k+1}-\bm{\varepsilon}^k\right),\left(\bm{\varepsilon}^{k+1}
+\bm{\varepsilon}^k\right)\right)\vspace{0.2 cm}\\\displaystyle
=\tau K_1\sum_{k=0}^{n-1}\left(\delta_{0,t}^{1-\alpha_1}\left(\mathscr{L}_y
\delta_x^2+\mathscr{L}_x\delta_y^2
\right)\bm{\varepsilon}^k,\left(\bm{\varepsilon}^{k+1}+\bm{\varepsilon}^k\right)\right)\vspace{0.2 cm}\\\displaystyle-
 \tau K_2\sum_{k=0}^{n-1}\left(\delta_{0,t}^{1-\alpha_2}\mathscr{L}_x\mathscr{L}_y \bm{\varepsilon}^k,
 \left(\bm{\varepsilon}^{k+1}+\bm{\varepsilon}^k\right)\right).
\end{array}\eqno(24)
$$
Obviously, the left hand side of equation (24) can be reduced to
$$
\begin{array}{ll}
\displaystyle
\sum_{k=0}^{n-1}\left(\mathscr{L}_x\mathscr{L}_y
\left(\bm{\varepsilon}^{k+1}-\bm{\varepsilon}^k\right),\left(\bm{\varepsilon}^{k+1}+\bm{\varepsilon}^k\right)\right)
=\left(\mathscr{L}_x\mathscr{L}_y\bm{\varepsilon}^{n},\bm{\varepsilon}^{n}\right)-
\left(\mathscr{L}_x\mathscr{L}_y\bm{\varepsilon}^{0},\bm{\varepsilon}^{0}\right).
\end{array}
$$
Applying Lemmas \ref{Lem.3.5}, \ref{Lem.4.5} and \ref{Lem.4.6} to the right hand side of equation (24) leads to
$$
\begin{array}{ll}
\displaystyle
\tau K_1\sum_{k=0}^{n-1}\left(\delta_{0,t}^{1-\alpha_1}\left(\mathscr{L}_y
\delta_x^2+\mathscr{L}_x\delta_y^2
\right)\bm{\varepsilon}^k,\left(\bm{\varepsilon}^{k+1}+\bm{\varepsilon}^k\right)\right)\vspace{0.2 cm}\\\displaystyle-
 \tau K_2\sum_{k=0}^{n-1}\left(\delta_{0,t}^{1-\alpha_2}\mathscr{L}_x\mathscr{L}_y
 \bm{\varepsilon}^k,\left(\bm{\varepsilon}^{k+1}+\bm{\varepsilon}^k\right)\right)\vspace{0.2 cm}\\\displaystyle
 =-\tau K_1\sum_{k=0}^{n-1}\left(\delta_{0,t}^{1-\alpha_1}\mathscr{Q}_{x+y}\bm{\varepsilon}^k,\mathscr{Q}_{x+y}
 \left(\bm{\varepsilon}^{k+1}+\bm{\varepsilon}^k\right)\right)\vspace{0.2 cm}\\\displaystyle-
 \tau K_2\sum_{k=0}^{n-1}\left(\delta_{0,t}^{1-\alpha_2}\mathscr{Q}_{xy}
 \bm{\varepsilon}^k,\mathscr{Q}_{xy}\left(\bm{\varepsilon}^{k+1}+\bm{\varepsilon}^k\right)\right)
 \vspace{0.2 cm}\\\displaystyle\leq0.
\end{array}
$$
Hence, combining the above two estimations gives
$$
\begin{array}{ll}
\displaystyle
\left(\mathscr{L}_x\mathscr{L}_y\bm{\varepsilon}^{n},\bm{\varepsilon}^{n}\right)\leq
\left(\mathscr{L}_x\mathscr{L}_y\bm{\varepsilon}^{0},\bm{\varepsilon}^{0}\right).
\end{array}
$$
Using Lemma \ref{Lem.4.4} again yields
$$
\begin{array}{ll}
\displaystyle
||\bm{\varepsilon}^{n}||\leq
\sqrt{3}||\bm{\varepsilon}^{0}||=\sqrt{3}||\bm{\rho}||.
\end{array}
$$
The proof is complete.
\end{proof}

Finally, we list the convergence result.
\begin{theorem}
Assume that the solution $u(x, y, t)$ of equation (17), together with (18) and (19) is suitably
smooth, and that $\{u_{i,j}^k|0\leq i\leq M_1,
0\leq j\leq M_2, 0 \leq k \leq N\}$ is the solution of the finite difference scheme
(21), together with (22) and (23). Let $e_{i,j}^k=u(x_i,y_j, t_k)-u_{i,j}^k$, then
$$
\begin{array}{lll}
\displaystyle
||\bm{e}^{n}||\leq 2C_2\exp\left(3T\right)\sqrt{6L_xL_yT}
\left(\tau^2+h_x^4+h_y^4\right).
\end{array}
$$
\end{theorem}

\begin{proof} The error equation reads as,
$$
\begin{array}{ll}
\displaystyle
\mathscr{L}_x\mathscr{L}_y
\left(e_{i,j}^{k+1}-e_{i,j}^k\right)
=\tau K_1\delta_{0,t}^{1-\alpha_1}\left(\mathscr{L}_y
\delta_x^2+\mathscr{L}_x\delta_y^2
\right)e_{i,j}^k- \tau K_2\delta_{0,t}^{1-\alpha_2}\mathscr{L}_x\mathscr{L}_y e_{i,j}^k\vspace{0.2 cm}\\
\displaystyle
 +\tau R_{i,j}^k,
\;\;\;\;\;\;1\leq i\leq M_1-1,\;1\leq j\leq M_2-1,\;0\leq k\leq N-1,\vspace{0.2 cm}\\
\displaystyle
u_{i,j}^k=0,\;\;\;\;\;1\leq i\leq M_1,\;1\leq j\leq M_2,\vspace{0.2 cm}\\
\displaystyle
u_{i,j}^k=\varphi\left(x_i,y_j,t_{k}\right),\;\;\;\;(x_i,y_j)\in\partial{\Omega}_{h_{x}h_{y}}.
\end{array}
$$
Simple calculations give
$$
\begin{array}{ll}
\displaystyle
\left(\mathscr{L}_x\mathscr{L}_y\bm{e}^{n},\bm{e}^{n}\right)\leq
\left(\mathscr{L}_x\mathscr{L}_y\bm{e}^{0},\bm{e}^{0}\right)+\tau\sum_{k=0}^{n-1}\left(\bm{R}^{k},\left(\bm{e}^{k+1}+\bm{e}^k\right)\right).
\end{array}
$$
Noticing
$$
\begin{array}{lll}
\displaystyle
\tau\sum_{k=0}^{n-1}\left(\bm{R}^{k},\left(\bm{e}^{k}+\bm{e}^k\right)\right)&\leq&\displaystyle\frac{\tau}{4}||\bm{e}^{n}||^2
+\frac{1}{2}\tau\sum_{k=0}^{n-1}
||\bm{e}^{k}||^2+
2\tau\sum_{k=0}^{n-1}
||\bm{R}^k||^2\vspace{0.2 cm}\\&\leq& \displaystyle
\frac{1}{4}||\bm{e}^{n}||^2
+\frac{1}{2}\tau\sum_{k=0}^{n-1}
||\bm{e}^{k}||^2+
2L_xL_yTC_2^2\left(\tau^2+h_x^4+h_y^4\right)^2,
\end{array}
$$
one has
$$
\begin{array}{lll}
\displaystyle
||\bm{e}^{n}||^2\leq
6\tau\sum_{k=0}^{n-1}
||\bm{e}^{k}||^2+
24L_xL_yTC_2^2\left(\tau^2+h_x^4+h_y^4\right)^2,
\end{array}
$$
i.e.,
$$
\begin{array}{lll}
\displaystyle
||\bm{e}^{n}||\leq
2C_2\exp\left(3T\right)\sqrt{6L_xL_yT}
\left(\tau^2+h_x^4+h_y^4\right).
\end{array}
$$
All this completes the proof.
\end{proof}

 \section{Numerical examples}

In this section, numerical experiments are carried out for the proposed numerical algorithms
 to illustrate our theoretical analysis.

\begin{example}
Take $u(t)=t^{2+\alpha}$, $t\in[0,1]$. Obviously, the exact value of $u(t)$ at $t=0.5$ is
$$
\begin{array}{lll}
\displaystyle \,_{\mathrm{RL}}{\mathrm{D}}_{0,t}^{\alpha}u\left(t\right)\left|\right._{t={0.5}}
=\frac{\Gamma(3+\alpha)}{\Gamma(2+2\alpha)}\left(\frac{1}{2}\right)^{1+2\alpha}.
\end{array}
$$
By using formula (4) to approximate function $u(t)$, the absolute error and convergence
 order are listed in Table \ref{tab.1}.
It can be seen from the table that the convergence order of formula (4) is almost 2,
which is in line with our theoretical analysis.

\begin{table}[!htbp]\renewcommand\arraystretch{1.2}
 \begin{center}
 \caption{ The absolute error and convergence order of Example 1 by numerical formula (4).}\label{tab.1}
 \vspace{0.5 cm}
 \begin{footnotesize}
\begin{tabular}{c c c c c c }\hline
  $\alpha$ &\; $\tau$\;&  \textrm{the absolute error}&\;\;\;\;$\textrm {the convergence order}$\\
  \hline
  $0.1 $& $\frac{1}{20}$ & 9.887230e-003 \;\; &  ---\\
  $$  & $\frac{1}{40}$&  3.031228e-003 \;\;&       1.7057 \\
  $$& $\frac{1}{80}$ &  8.357173e-004 \;\; &      	1.8588 \\
  $$& $\frac{1}{160}$ &  2.184140e-004	 \;\; &     1.9359\\
   $$& $\frac{1}{320}$ & 5.576440e-005 \;\; &    1.9696\\
 \hline
   $0.3 $& $\frac{1}{20}$ & 8.294679e-003 \;\;& ---\\
  $$  & $\frac{1}{40}$& 2.177047e-003	 \;\;  &   	1.9298\\
   $$  & $\frac{1}{80}$&5.554318e-004  \;\;  &   	1.9707\\
  $$& $\frac{1}{160}$ &1.401845e-004 \;\; &      	1.9863\\
   $$& $\frac{1}{320}$ & 3.520843e-005 \;\; &    	1.9933\\
 \hline
  $0.5 $& $\frac{1}{20}$ &7.611138e-003\;\; &     ---\\
  $$  & $\frac{1}{40}$& 1.903986e-003   \;\; &     1.9991\\
  $$  & $\frac{1}{80}$&4.759995e-004    \;\; &       2.0000\\
  $$& $\frac{1}{160}$ &  1.190001e-004\;\;  &  2.0000  \\
   $$& $\frac{1}{320}$ & 	2.975004e-005	 \;\; &    2.0000\\
 \hline
   $0.7 $& $\frac{1}{20}$ &5.717212e-003 \;\; &      ---\\
  $$  & $\frac{1}{40}$& 1.389258e-003 \;\; &      2.0410 \\
   $$  & $\frac{1}{80}$& 3.429045e-004\;\; &     2.0184 \\
  $$& $\frac{1}{160}$ & 8.520286e-005   \;\; &   2.0088\\
   $$& $\frac{1}{320}$ &  2.123685e-005 \;\; &   2.0043\\
 \hline
$0.9 $& $\frac{1}{20}$ &	2.151549e-003 \;\; &      ---\\
  $$  & $\frac{1}{40}$& 5.116059e-004	 \;\; &     	2.0723 \\
   $$  & $\frac{1}{80}$&  1.249915e-004   \;\; &  2.0332\\
  $$& $\frac{1}{160}$ & 3.090217e-005	    \;\; &  2.0160\\
   $$& $\frac{1}{320}$ &  7.683366e-006 \;\; &    2.0079\\
 \hline
\end{tabular}
 \end{footnotesize}
 \end{center}
 \end{table}
\end{example}

\begin{example}
Consider the following one-dimensional fractional cable equation:
$$
\begin{array}{rrr}
\displaystyle
\frac{\partial u(x,t)}{\partial t}=\frac{1}{\pi^8}\,_{\mathrm{RL}}{{{\mathrm{D}}}}_{0,t}^{1-\alpha_1}
\frac{\partial^2u(x,t) }{\partial x^2}
- \,_{\mathrm{RL}}{{{\mathrm{D}}}}_{0,t}^{1-\alpha_2} u(x,t) +f(x,t),\;(x,t)\in(0,1]\times(0,1],
\end{array}
$$
with initial value condition
$$
\begin{array}{lll}
\displaystyle
u(x,0)=0,\;\;x\in[0,1],
\end{array}
$$
and boundary value conditions
$$
\begin{array}{lll}
\displaystyle
u(0,t)=0,\;u(1,t)=0,\;\;t\in(0,1],
\end{array}
$$
where the source term is $f(x,t)=2\left(t+\frac{t^{1+\gamma_1}}{\pi^6\Gamma(2+\gamma_1)}+\frac{t^{1+\gamma_2}}{\Gamma(2+\gamma_2)}\right)\sin\pi x$.
The analytical solution of the above system is $u(x, t) = t^2 \sin\pi x$.
\end{example}

\begin{figure}[!htbp]
\begin{center}
 \includegraphics[width=10 cm]{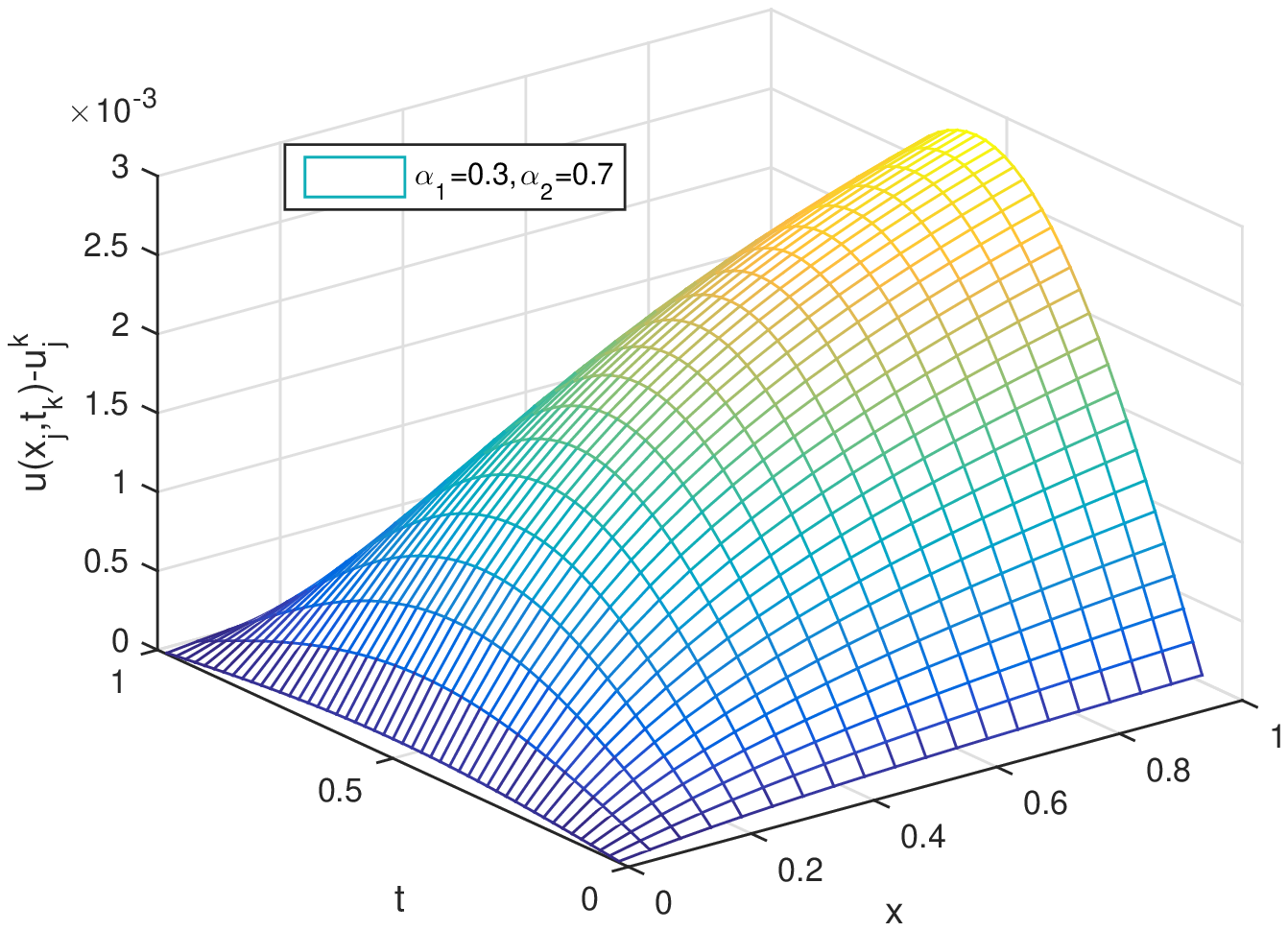}\\
  \caption{The evolutions of the absolute errors obtained by the numerical algorithm (11), together with (12) and (13) for
   $(\alpha_1,\alpha_2)=(0.3,0.7)$ when $\tau=\frac{1}{20}$ and $h=\frac{1}{50}$.
  }\label{Fig.3}
  \end{center}
\end{figure}

\begin{figure}[!htbp]
\begin{center}
 \includegraphics[width=10 cm]{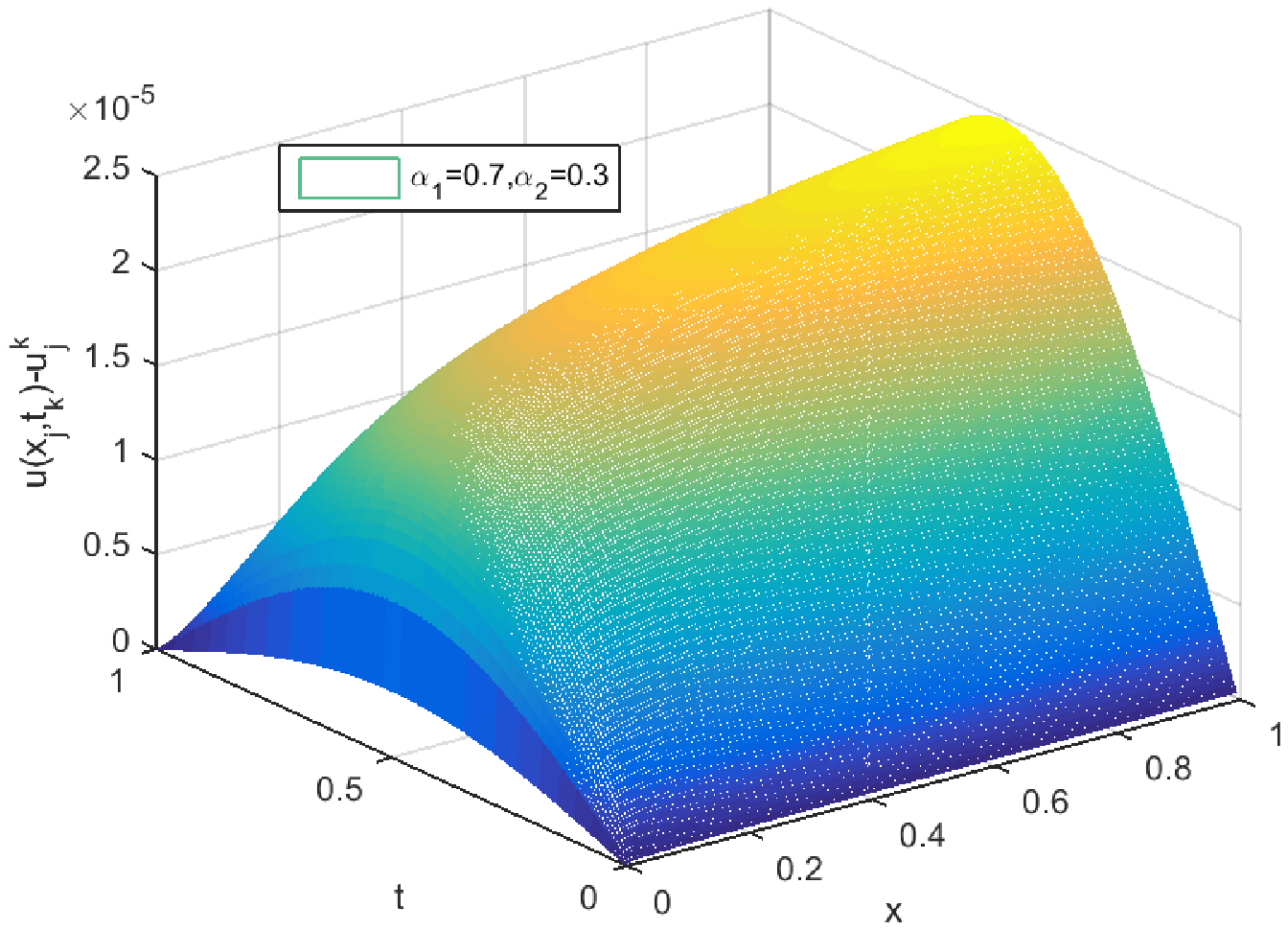}\\
  \caption{The evolutions of the absolute errors obtained by the numerical algorithm (11), together with (12) and (13) for
   $(\alpha_1,\alpha_2)=(0.7,0.3)$ when $\tau=\frac{1}{200}$ and $h=\frac{1}{200}$.
  }\label{Fig.4}
  \end{center}
\end{figure}

In this numerical test, we present the absolute errors and the corresponding temporal and spatial
convergence orders in Table \ref{tab.2} for different $\alpha_1$ and $\alpha_2$ by using finite difference scheme (11), together with (12) and (13),
which verifies that the second-order accuracy in time
and fourth-order accuracy in space direction are obtained. Meanwhile,
the evolutions of the absolute errors were depicted in Figs. \ref{Fig.3} and \ref{Fig.4} for different orders $\alpha_1,\alpha_2$
 and stepsizes $\tau,h$.
Obviously,
all of the above numerical results are in accordance with our stability and convergence analysis of the proposed numerical algorithm (11), together with (12) and (13).

\begin{table}\renewcommand\arraystretch{1.2}
 \begin{center}
 \caption{ The absolute errors (TAE), temporal convergence
order (TCO) and spatial convergence order (SCO) of Example 2 by using difference scheme (11), together with (12) and (13).}\label{tab.2}
 \vspace{0.3 cm}
 \begin{footnotesize}
\begin{tabular}{c c c c c c }\hline
   & \;\;\;\; \;\; &\;\;\;\;\;\;&\;\;\;\;\;\;&\;\;\;\;\;\;&\;\;\;\;\;\;\\
  $(\alpha_1,\alpha_2)$ &\;$\tau$,\;$h$\;&  \;\;TAE\;\;\;\; &\;\;\;\; TCO\;\;\;\; & SCO\\\hline \vspace{0.1 cm}
  $(0.2,0.8) $& $\tau=\frac{1}{5},h=\frac{1}{5}$ &  4.017482e-02		&  ---&  ---\\ \vspace{0.1 cm}
  $$  & $\tau=\frac{1}{20},h=\frac{1}{10}$&       3.443878e-03	 &  1.7721&   3.5442\\ \vspace{0.1 cm}
  $$& $\tau=\frac{1}{45},h=\frac{1}{15}$ &      7.010198e-04	&   1.9630&   3.9259\\ \vspace{0.1 cm}
$$  & $\tau=\frac{1}{80},h=\frac{1}{20}$&        2.254316e-04	 &    1.9718&  3.9437\\ \vspace{0.1 cm}
  $$& $\tau=\frac{1}{125},h=\frac{1}{25}$ &       9.258968e-05	&  1.9939 &   3.9877\\ \hline\vspace{0.1 cm}
   $(0.4,0.6) $& $\tau=\frac{1}{5},h=\frac{1}{5}$ &   3.636258e-02	&  ---&  ---\\ \vspace{0.1 cm}
  $$  & $\tau=\frac{1}{20},h=\frac{1}{10}$&         2.709682e-03 &   1.8731&   3.7463\\ \vspace{0.1 cm}
  $$& $\tau=\frac{1}{45},h=\frac{1}{15}$ &       5.399822e-04	&    1.9891&   3.9783\\ \vspace{0.1 cm}
$$  & $\tau=\frac{1}{80},h=\frac{1}{20}$&        1.725672e-04	 &    1.9827&   3.9653\\ \vspace{0.1 cm}
  $$& $\tau=\frac{1}{125},h=\frac{1}{25}$ &       7.068420e-05	&  2.0000 &   4.0000\\ \hline\vspace{0.1 cm}
  $(0.5,0.5) $& $\tau=\frac{1}{5},h=\frac{1}{5}$ &  3.552136e-02	&  ---&  ---\\ \vspace{0.1 cm}
  $$  & $\tau=\frac{1}{20},h=\frac{1}{10}$&        2.583580e-03	 &   1.8906&  3.7812\\ \vspace{0.1 cm}
  $$& $\tau=\frac{1}{45},h=\frac{1}{15}$ &      5.125456e-04&    1.9947&   3.9893\\ \vspace{0.1 cm}
$$  & $\tau=\frac{1}{80},h=\frac{1}{20}$&        1.635583e-04	 &  1.9852&   3.9704\\ \vspace{0.1 cm}
  $$& $\tau=\frac{1}{125},h=\frac{1}{25}$ &       6.694813e-05	& 2.0015 &   4.0030\\ \hline\vspace{0.1 cm}
   $(0.6,0.4) $& $\tau=\frac{1}{5},h=\frac{1}{5}$ &   3.463007e-02	&  ---&  ---\\ \vspace{0.1 cm}
  $$  & $\tau=\frac{1}{20},h=\frac{1}{10}$&       2.513471e-03 &  1.8921&   3.7843\\ \vspace{0.1 cm}
  $$& $\tau=\frac{1}{45},h=\frac{1}{15}$ &       4.969394e-04	&   1.9989&   3.9978\\ \vspace{0.1 cm}
$$  & $\tau=\frac{1}{80},h=\frac{1}{20}$&       1.583867e-04 &    1.9873&  3.9746\\ \vspace{0.1 cm}
  $$& $\tau=\frac{1}{125},h=\frac{1}{25}$ &       6.479215e-05&  2.0029 &   4.0057\\ \hline\vspace{0.1 cm}
  $(0.8,0.2) $& $\tau=\frac{1}{5},h=\frac{1}{5}$ &   4.961263e-02	&  ---&  ---\\ \vspace{0.1 cm}
  $$  & $\tau=\frac{1}{20},h=\frac{1}{10}$&        2.863711e-03	 &2.0574   & 4.2960  \\ \vspace{0.1 cm}
  $$& $\tau=\frac{1}{45},h=\frac{1}{15}$ &      5.157444e-04	& 2.1139   & 4.0857  \\ \vspace{0.1 cm}
$$  & $\tau=\frac{1}{80},h=\frac{1}{20}$&       1.542104e-04	 &  2.0983   &  4.0304 \\ \vspace{0.1 cm}
  $$& $\tau=\frac{1}{125},h=\frac{1}{25}$ &      6.299963e-05	& 2.0059 &  4.0127  \\ \hline\vspace{0.1 cm}
\end{tabular}
 \end{footnotesize}
 \end{center}
 \end{table}

\begin{example}
Consider the following two-dimensional fractional cable equation:
$$
\begin{array}{rrr}
\displaystyle
\frac{\partial u(x,y,t)}{\partial t}=\frac{1}{\pi^8}\,_{\mathrm{RL}}{{{\mathrm{D}}}}_{0,t}^{1-\alpha_1}
\left(\frac{\partial^2u(x,y,t) }{\partial x^2}+\frac{\partial^2u(x,y,t) }{\partial y^2}\right)
- \,_{\mathrm{RL}}{{{\mathrm{D}}}}_{0,t}^{1-\alpha_2} u(x,y,t)\vspace{0.2 cm}\\ +f(x,y,t),\;(x,y;t)\in\Omega\times[0,1],
\end{array}
$$
with initial value condition
$$
\begin{array}{lll}
\displaystyle
u(x,y,0)=0,\;\;(x,y)\in\Omega,
\end{array}
$$
and boundary value conditions
$$
\begin{array}{lll}
\displaystyle
u(x,y,t)=0,\;(x,y)\in\partial\Omega,\;t\in(0,1],
\end{array}
$$
where the spatial domain $\Omega=(0,1)\times(0,1)$, and the source term is chosen as $f(x,y,t)=2\left(t+\frac{2t^{1+\gamma_1}}{\pi^6\Gamma(2+\gamma_1)}+\frac{t^{1+\gamma_2}}{\Gamma(2+\gamma_2)}\right)\sin\pi x\sin\pi y$.
The analytical solution is $u(x, t) = t^2 \sin\pi x\sin\pi y$.
\end{example}

\begin{figure}[!htbp]
\begin{center}
 \includegraphics[width=10 cm]{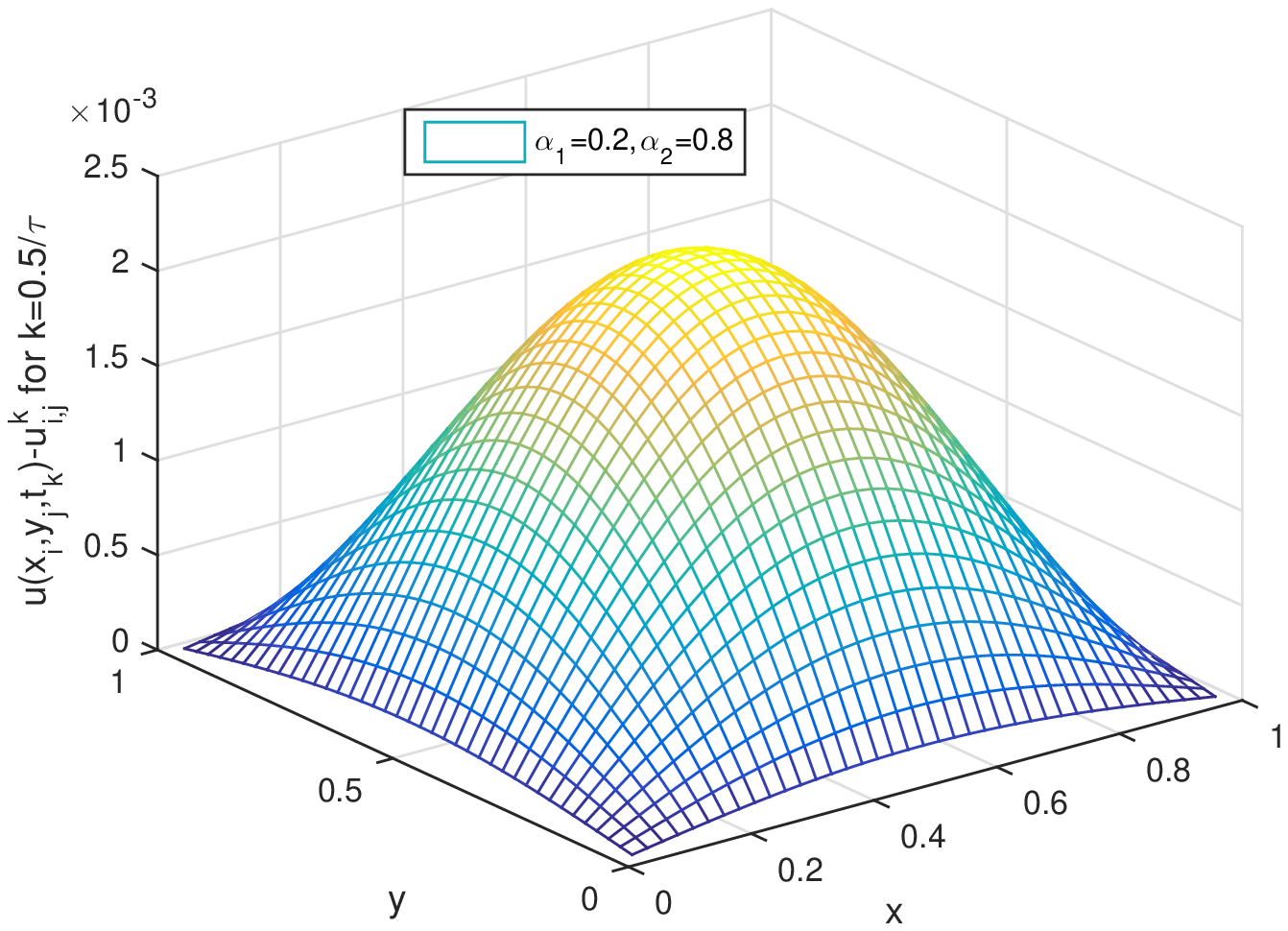}\\
  \caption{The evolutions of the absolute errors obtained by the numerical algorithm (21), together with (22) and (23) for
   $(\alpha_1,\alpha_2)=(0.2,0.8)$ at $t=0.5$ when $\tau=\frac{1}{20}$ and $h=\frac{1}{40}$.
  }\label{Fig.5}
  \end{center}
\end{figure}

\begin{figure}[!htbp]
\begin{center}
 \includegraphics[width=10 cm]{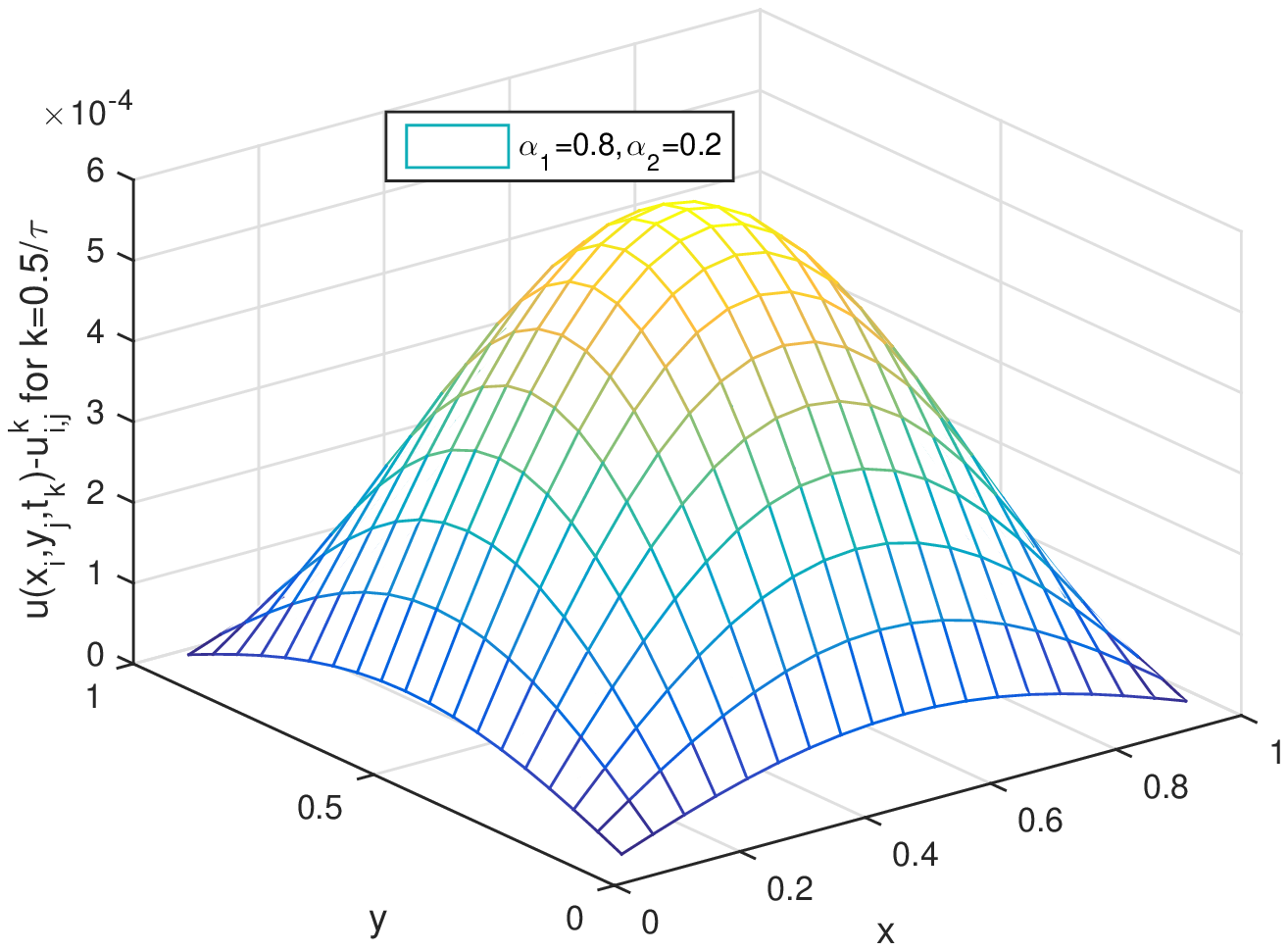}\\
  \caption{The evolutions of the absolute errors obtained by the numerical algorithm (21), together with (22) and (23) for
   $(\alpha_1,\alpha_2)=(0.8,0.2)$ at $t=0.5$ when $\tau=\frac{1}{40}$ and $h=\frac{1}{20}$.
  }\label{Fig.6}
  \end{center}
\end{figure}

 Table \ref{tab.3} lists the computed errors,
 the temporal and spatial convergence orders respectively, which shows that the convergence order
  of our scheme (21)--(23) is $\mathcal{O}\left(\tau^2+h_x^4+h_y^4\right)$. Finally, Figs. \ref{Fig.5} and \ref{Fig.6} show the evolutions of the absolute errors for
   Example 5.3 with different orders, temporal and spatial stepsizes, respectively.
It can be seen that the numerical results are in good
agreement with the theoretical results.

\begin{table}\renewcommand\arraystretch{1.2}
 \begin{center}
 \caption{ The absolute errors (TAE), temporal convergence
order (TCO) and spatial convergence order (SCO) of Example 3 by using difference scheme (21), together with (22) and (23).}\label{tab.3}
 \vspace{0.3 cm}
 \begin{footnotesize}
\begin{tabular}{c c c c c c }\hline
   & \;\;\;\; \;\; &\;\;\;\;\;\;&\;\;\;\;\;\;&\;\;\;\;\;\;&\;\;\;\;\;\;\\
  $(\alpha_1,\alpha_2)$ &\;$\tau$,\;$h_x,h_y$\;&  \;\;TAE\;\;\;\; &\;\;\;\; TCO\;\;\;\; & SCO\\\hline \vspace{0.1 cm}
  $(0.2,0.8) $& $\tau=\frac{1}{5},h_x=h_y=\frac{1}{5}$ &  3.822113e-02			&  ---&  ---\\ \vspace{0.1 cm}
  $$  & $\tau=\frac{1}{20},h_x=h_y=\frac{1}{10}$&       3.444335e-03		 &  1.7360&   3.4721\\ \vspace{0.1 cm}
  $$& $\tau=\frac{1}{45},h_x=h_y=\frac{1}{15}$ &      6.972544e-04		&   1.9698&   3.9395\\ \vspace{0.1 cm}
$$  & $\tau=\frac{1}{80},h_x=h_y=\frac{1}{20}$&        2.254541e-04		 &    1.9623&  3.9246\\ \vspace{0.1 cm}
  $$& $\tau=\frac{1}{125},h_x=h_y=\frac{1}{25}$ &       2.254541e-04		& 1.9983 &  3.9966\\ \hline\vspace{0.1 cm}
   $(0.4,0.6) $& $\tau=\frac{1}{5},h_x=h_y=\frac{1}{5}$ &  3.460010e-02	&  ---&  ---\\ \vspace{0.1 cm}
  $$  & $\tau=\frac{1}{20},h_x=h_y=\frac{1}{10}$&        2.710859e-03	&   1.8370&  3.6740\\ \vspace{0.1 cm}
  $$& $\tau=\frac{1}{45},h_x=h_y=\frac{1}{15}$ &      5.372495e-04		&   1.9959&  3.9919\\ \vspace{0.1 cm}
$$  & $\tau=\frac{1}{80},h_x=h_y=\frac{1}{20}$&       1.726387e-04		 &   1.9731&   3.9462\\ \vspace{0.1 cm}
  $$& $\tau=\frac{1}{125},h_x=h_y=\frac{1}{25}$ &      7.057381e-05	&  2.0044 &   4.0088\\ \hline\vspace{0.1 cm}
  $(0.5,0.5) $& $\tau=\frac{1}{5},h_x=h_y=\frac{1}{5}$ & 3.379936e-02	&  ---&  ---\\ \vspace{0.1 cm}
  $$  & $\tau=\frac{1}{20},h_x=h_y=\frac{1}{10}$&        2.584901e-03		 &  1.8544&  3.7088\\ \vspace{0.1 cm}
  $$& $\tau=\frac{1}{45},h_x=h_y=\frac{1}{15}$ &      5.099957e-04	&   2.0015& 4.0029\\ \vspace{0.1 cm}
$$  & $\tau=\frac{1}{80},h_x=h_y=\frac{1}{20}$&       1.636407e-04		 & 1.9757&  3.9513\\ \vspace{0.1 cm}
  $$& $\tau=\frac{1}{125},h_x=h_y=\frac{1}{25}$ &    6.684963e-05			& 2.0059 &  4.0119\\ \hline\vspace{0.1 cm}
   $(0.6,0.4) $& $\tau=\frac{1}{5},h_x=h_y=\frac{1}{5}$ & 3.295086e-02			&  ---&  ---\\ \vspace{0.1 cm}
  $$  & $\tau=\frac{1}{20},h_x=h_y=\frac{1}{10}$&      2.514903e-03	 & 1.8559&   3.7117\\ \vspace{0.1 cm}
  $$& $\tau=\frac{1}{45},h_x=h_y=\frac{1}{15}$ &     4.945015e-04		&  2.0056&   4.0113\\ \vspace{0.1 cm}
$$  & $\tau=\frac{1}{80},h_x=h_y=\frac{1}{20}$&     1.584782e-04	 &  1.9778& 3.9555\\ \vspace{0.1 cm}
  $$& $\tau=\frac{1}{125},h_x=h_y=\frac{1}{25}$ &    6.470169e-05	&  2.0073 &   4.0146\\ \hline\vspace{0.1 cm}
  $(0.8,0.2) $& $\tau=\frac{1}{5},h_x=h_y=\frac{1}{5}$ &  4.721574e-02		&  ---&  ---\\ \vspace{0.1 cm}
  $$  & $\tau=\frac{1}{20},h_x=h_y=\frac{1}{10}$&      2.866291e-03		 &   2.0210&   4.0420\\ \vspace{0.1 cm}
  $$& $\tau=\frac{1}{45},h_x=h_y=\frac{1}{15}$ &    5.133508e-04			&   2.1208&  4.2416\\ \vspace{0.1 cm}
$$  & $\tau=\frac{1}{80},h_x=h_y=\frac{1}{20}$&    1.543393e-04			 &    2.0888&   4.1775\\ \vspace{0.1 cm}
  $$& $\tau=\frac{1}{125},h_x=h_y=\frac{1}{25}$ &    6.292841e-05		&  2.0103 &  4.0205\\ \hline\vspace{0.1 cm}
\end{tabular}
 \end{footnotesize}
 \end{center}
 \end{table}

\section{Conclusions}

In this paper, two classed high-order numerical algorithms are derived to solve
 the one- and two-dimensional fractional Cable equations
based on the derived new second-order difference operator in time direction and the
compact techniques in space direction. By using the energy method, we
proved that our difference schemes are unconditionally stable to the initial values.
 The temporal, spatial convergence orders can reach two and four respectively.
 Finally, numerical examples are given
to show the effectiveness of the derived numerical algorithms.


\begin{thebibliography}{}

\bibitem{A} A.A. Alikhanov, A new difference scheme for the fractional diffusion equation, J. Comput. Phys., 280
(2015), 424--438.
\bibitem{BDA}D. Baleanu, O. Defterli, O.P. Agrawal, A central difference numerical scheme for fractional
optimal control problems, J. Vib. Control., 15 (2009), 583--597.

\bibitem{BWM} D. Benson, S.W. Wheatcraft, M.M. Meerschaert, The fractional-order governing equation of L\'{e}vy motion,
 Water Resour. Res., 36 (2000), 1413--1423.

\bibitem{C}M. Cui, Convergence analysis of high-order compact alternating direction implicit schemes for the
 two-dimensional time fractional diffusion equation, Numer. Algorithms, 62 (2013), 383--409.

\bibitem{DL1}  H.F. Ding,  C.P. Li, Y.Q. Chen, High-order Algorithms
for Riesz Derivative and Their Applications (I), Abstract and Applied
Analysis, Article ID 653797, 2014, 17 pages.

\bibitem{DL2}
H.F. Ding,  C.P. Li, Y.Q. Chen, High-order
algorithms for Riesz derivative and their applications
(II), J. Comput. Phys., 293  (2015), 218--237.

\bibitem{ER} V.J. Ervin, J.P. Roop, Variational formulation for the stationary fractional advection dispersion
 equation,
Numer. Meth. Part. D. E.,
22 (2005), 558--576.

\bibitem{GS}G.H. Gao, Z.Z. Sun, A compact finite difference scheme for the fractional
 sub-diffusion equations, J. Comput. Phys., 230 (2011), 586--595.

\bibitem{HZ}X.L. Hu, L.M. Zhang,
Implicit compact difference schemes for the fractional cable equation, Appl. Math. Model.,
36 (2012), 4027--4043.

\bibitem{J}X. Jin, Preconditioning Techniques for Toeplitz Systems, Higher Education Press, Beijing, 2010.

\bibitem{JLPZ} B. Jin, R. Lazarov, J. Pasciak, Z. Zhou, Error analysis of a finite element method for
the space-fractional parabolic equation, SIAM. J. Numer. Anal., 52 (2014), 2272--2294.

\bibitem{JS}C.C. Ji, Z.Z. Sun, The high-order compact numerical algorithms for
 the two-dimensional fractional sub-diffusion equation, Appl. Math. Comput., 269 (2015), 775--791.

\bibitem{L}A.J. Laub, Matrix Analysis for Scientists and Engineers, Society for Industrial
 and Applied Mathematics, Philadelphia, PA, 2005.

\bibitem{LAT} F. Liu, V. Anh, I. Turner, Numerical solution of the space fractional
 Fokker-Planck equation, J. Comput. Appl. Math., 166 (2004), 209--219.

\bibitem{LD}C.P. Li, H.F. Ding, Higher order finite difference method for the reaction
 and anomalous-diffusion
equation, Appl. Math. Model., 38 (2014), 3802--3821.

\bibitem{LHW}
T.A.M. Langlans, B. Henry, S. Wearne, Solution of a fractional Cable equation: finite case,
 Preprint, Submitted to Elsevier Science http://www.maths.unsw.edu.au/applied/filed/2005/amr05-33 (2005).

\bibitem{LLX}
Y.M. Lin, X.J. Li, C.J. Xu, Finite difference/spectral approximations for the fractional Cable equation,
Math. Comput., 80 (2011), 1369--1396 .

\bibitem{MK1} R. Metler, J. Klafter, The random walk's guide to anomalous diffusion: a fractional
 dynamics approach, Phys. Rep., 339 (2000), 1--77.

\bibitem{MK2} R. Metler, J. Klafter, The restaurant at the end of random walk: recent developments
 in the description of anomalous transport by fractional dynamics,
J. Phys. A., 37 (2004), R161--R208.

\bibitem{MRGS} F. Mainardi, M. Raberto, R. Gorenflo, E. Scalas, Fractional calculus and
 continuous-time finance I: the waiting-time distribution, Physica A., 287 (2000),
468--481.

\bibitem{MT} M.M. Meerschaert, C. Tadjeran, Finite difference approximations for fractional
 advection-dispersion flow equations, J. Comput. Appl. Math., 172 (2004),
65--77.

\bibitem{MY}J. Quintana-Murillo, S.B. Yuste, An explicit numerical method for the fractional cable equation, Int. J.
Differ. Equ., Article ID 231920, 2011, 12 pages.

\bibitem{P}
I. Podlubny, Fractional Differential Equations, Academic Press, San Diego, 1999.

\bibitem{PCSCJ}I. Podlubny, A. Chechkin, T. Skovranek, Y. Chen, B. M. V. Jara, Matrix approach to discrete fractional calculus.
 II. Partial fractional differential equations, J. Comput. Phys., 228 (2009), 3137--3153.

\bibitem{QSS} A. Quarteroni, R. Sacco, F. Saleri, Numerical
Mathematics, Springer, New York, 2007.

\bibitem{SK} I.M. Sokolov, J. Klafter, From diffusion to anomalous diffusion: a century after Einstein's Brownian motion,
 Chaos, 15 (2005), 026103.

\bibitem{SKM}
 S.G. Samko, A.A. Kilbas, O.I. Marichev, Fractional Integrals and Derivatives: Theory and Applications,
Gordon and Breach, 1993.

\bibitem{SLLA}S. Shen, F. Liu, Q. Liu, V. Anh, Numerical simulation of anomalous
 infiltration in porous media, Numer. Algorithms, 68 (2015), 443--454.

\bibitem{WG}
T. Wang, B. Guo, A robust semi-explicit difference scheme for the Kuramoto-Tsuzuki equation,
 J. Comput. Appl. Math., 233 (2009), 878--888.

\bibitem{WV}Z.B. Wang, S. Vong, Compact difference schemes for the modified anomalous fractional
sub-diffusion equation and the fractional diffusion-wave equation, J. Comput. Phys., 277 (2014), 1--15.

\bibitem{WY}H. Wang, D. Yang, Well posedness of variable-coefficient conservative fractional elliptic
 differential equations, SIAM. J. Numer. Anal., 51 (2013), 1088--1107.

\bibitem{ZLTA} M. Zheng, F. Liu, I. Turner, V. Anh, A novel high order space-time spectral method for the
 time-fractional Fokker-Planck equation, SIAM. J. Sci. Comput., 37 (2015), 701--724.

\bibitem{ZLTA}P. Zhuang, F. Liu,  I. Turner, V. Anh,
Galerkin finite element method and error analysis
for the fractional Cable equation, Numer. Algorithms, DOI 10.1007/s11075-015-0055-x.

\bibitem{ZK} M. Zayernouri, G.E. Karniadakis, Discontinuous spectral element methods for
time-and space-fractional advection equations, SIAM. J. Sci. Comput., 36 (2014), B684--B707.
\end{thebibliography}
  \end{document}